\documentclass[a4paper,12pt,nosumlimits]{amsart}
\usepackage{latexsym,amsmath,amssymb,color,mathrsfs,enumerate,bbm}
\usepackage[english]{babel}

\numberwithin{equation}{section}

\title[Differential Geometry of the Quantum Hopf Fibration]{Differential and Twistor Geometry of the \\ Quantum Hopf Fibration}
\author{Simon Brain} \address{Unit\'{e} de Recherche en Math\'{e}matiques, Universit\'{e} du Luxembourg (Campus Kirchberg), 6 rue Richard Coudenhove--Kalergi, L-1359 Luxembourg, Grand Duchy of Luxembourg}\email{simon.brain@uni.lu}

\author{Giovanni Landi}
\address{Dipartimento di Matematica, 
Universit\`{a} di Trieste, Via A. Valerio 12/1, 34127
Trieste, Italy and INFN, Sezione di Trieste, Trieste,
Italy}\email{landi@units.it}

\date{v2: 7 September 2012}

\newtheorem{thm}{Theorem}[section]

\newtheorem{lem}[thm]{Lemma}
\newtheorem{prop}[thm]{Proposition}
\newtheorem{rem}[thm]{Remark}

\theoremstyle{definition}
\newtheorem{defn}[thm]{Definition}
\newtheorem{example}[thm]{Example}

\newcommand{\id}{\textup{id}}
\newcommand{\D}{\textup{d}}
\newcommand{\iq}{q^{-1}}
\newcommand{\la}{\langle}
\newcommand{\ra}{\rangle}

\newcommand{\n}{\nabla}
\newcommand{\ep}{\epsilon}
\newcommand{\M}{\textup{M}}
\newcommand{\bfx}{\mathbf{x}}
\newcommand{\bfa}{\mathbf{a}}
\newcommand{\bfb}{\mathbf{b}}
\newcommand{\bfc}{\mathbf{c}}
\newcommand{\bfd}{\mathbf{d}}

\newcommand{\nn}{\nonumber}
\newcommand{\ASU}{\A[\SU_q(2)]}
\newcommand{\hV}{\widehat V}
\newcommand{\tV}{\widetilde V}
\newcommand{\hH}{\widehat H}
\newcommand{\tH}{\widetilde H}
\newcommand{\tM}{\widetilde M}
\newcommand{\hM}{\widehat M}

\newcommand{\cop}{\textup{cop}}

\newcommand{\pp}{{\sf p}}
\newcommand{\qp}{{\sf q}}
\newcommand{\Pp}{{\sf P}}

\newcommand{\sx}{\mathrm{x}}

\newcommand{\sv}{\mathsf{v}}
\newcommand{\su}{\sf u}
\newcommand{\sfu}{\sf u}

\newcommand{\sa}{\mathrm{a}}
\renewcommand{\sb}{\mathrm{b}}
\renewcommand{\sc}{\mathrm{c}}
\newcommand{\sd}{\mathrm{d}}

\newcommand{\st}{\mathsf{t}}

\newcommand{\A}{\mathcal{A}}
\newcommand{\E}{\mathcal{E}}

\newcommand{\Sp}{\textup{Sp}}

\newcommand{\U}{\textup{U}}
\newcommand{\SU}{\textup{SU}}

\newcommand{\CP}{\mathbb{C} \mathbb{P}}

\newcommand{\HH}{\mathbb{H}}
\newcommand{\ZZ}{\mathbb{Z}}
\newcommand{\C}{\mathbb{C}}
\newcommand{\RR}{\mathbb{R}}

\renewcommand{\o}{{}_{\scriptscriptstyle(1)}}
\renewcommand{\t}{{}_{\scriptscriptstyle(2)}}
\newcommand{\thr}{{}_{\scriptscriptstyle(3)}}

\newcommand{\uo}{{}^{{\scriptscriptstyle(1)}}}

\newcommand{\bz}{{}^{{\scriptscriptstyle(0)}}}

\def\rbiprod{{\cdot\kern-.33em\triangleright\!\!\!<}}
\def\lbiprod{{>\!\!\!\triangleleft\kern-.33em\cdot}}

\newcommand{\tDelta}{\widetilde{\Delta}}
\newcommand{\tS}{\widetilde{S}}

\begin{document}

\begin{abstract}
We study a quantum version of the $\textrm{SU}(2)$ Hopf fibration $S^7
\to S^4$ and its associated twistor geometry. Our quantum sphere
$S^7_q$ arises as the unit sphere inside a $q$-deformed quaternion
space $\mathbb{H}^2_q$. The resulting four-sphere $S^4_q$ is a quantum
analogue of the quaternionic projective space $\mathbb{HP}^1$. 
The quantum fibration is endowed with compatible non-universal differential calculi. By investigating the quantum symmetries of the fibration, we obtain the geometry of the corresponding twistor space $\mathbb{CP}^3_q$ and use it to study a system of anti-self-duality equations on $S^4_q$, for which we find an `instanton' solution coming from the natural projection defining the tautological bundle over $S^4_q$.
\end{abstract}

\maketitle
\tableofcontents

\linespread{1.1} 
\parskip 1ex

\section{Introduction}

The Hopf fibration over the standard quantum two-sphere $S^2_q$ was constructed well over a decade ago \cite{brmaj} and provided the very first interesting example of a quantum principal bundle with non-trivial differential structures. The underlying geometry of that particular example is a $q$-deformation of the canonical line bundle over the one-dimensional complex projective space $\CP^1$. In the present paper we introduce and study what we deem to be a $q$-deformation of the analogous bundle over the quaternionic projective space $\mathbb{H}\mathbb{P}^1$.

In classical geometry, the quaternionic Hopf fibration is nothing other than the principal $\SU(2)$ bundle $S^7\to S^4$. Recent years have witnessed a variety of attempts to generalise the geometry of this fibration to the framework of quantum group theory. The setting of `isospectral deformations', in which the deformation parameter is a complex number of modulus one \cite{cl}, has enjoyed a degree of success in this direction \cite{lvs:pfns}. On the other hand, the quest to find such examples in the setting of $q$-deformed geometry, where $q\ne 0$ is a real deformation parameter, has proven to be much more difficult. Interesting examples of $q$-deformed Hopf bundles over quantum four-spheres appeared in \cite{bct,lpr:qsb}, although both were described only at the topological level and their finite-dimensional differential structures have so far proven elusive.

In the present paper we start from observing that, if one seeks to find a $q$-deformation of the quaternionic Hopf fibration, then one should begin with a deformation of the quaternions themselves. Our route goes {\em via} the quantum Euclidean space $\RR^4_q$ introduced in \cite{ma:twist}. 
We endow this noncommutative space $\RR^4_q$ with a quaternionic structure, from which we construct a quantum seven-sphere $S^7_q$ carrying an action of the quantum group $\SU_q(2)$ whose quotient is a certain quantum four-sphere $S^4_q$. 

It is precisely the presence of the quaternionic structure which allows us to find finite-dimensional differential structures on our quantum spheres and to construct a noncommutative `instanton gauge configuration', a quaternionic analogue of the Dirac $q$-monopole found in \cite{brmaj}, which we interpret as a deformation of the BPST instanton of \cite{BPST}. Indeed, the secondary stream of the present paper is to describe the differential geometry of the quantum twistor fibration $\CP^3_q\to S^4_q$. In the classical case \cite{ma:gymf}, this twistor fibration encodes the very nature of the anti-self-dual Yang-Mills equations on the Euclidean sphere $S^4$. We find here that in the $q$-deformed setting, this is indeed the case for the sphere 
$S^4_q$, which may also be characterised as being Euclidean. The quantum spheres and projective space constructed herein appear to be different from those studied in \cite{vaks,welk2} or the more recent examples in \cite{dlm,lpr:qsb}, although there are several similarities as pointed out below.

The paper is organised as follows. In \S\ref{se:prelims} we review the basic concepts that we shall need regarding the structure of quantum principal bundles. In particular, we recall the quantum group $\SU_q(2)$ and its representation theory, together with its four-dimensional bicovariant differential calculus. In \S\ref{se:hopf-fib} we present the quantum Hopf fibration itself: a quantum principal bundle $S^7_q\to S^4_q$ with structure quantum group $\SU_q(2)$. We also find the noncommutative twistor space $\CP^3_q$ and the associated fibration $\CP^3_q\to S^4_q$.

In \S\ref{se:symm} we compute the group of quantum symmetries of our Hopf fibration, which we use to realise the noncommutative spaces $S^7_q$ and $\CP^3_q$ as quantum homogeneous spaces thereof. This is where our approach differs from those of \cite{bct,lpr:qsb}, both of which assume from the outset a quantum group of symmetries and then look to derive from this an appropriate fibration of homogeneous spaces. Nevertheless, just as in \cite{lpr:qsb}, our noncommutative four-sphere $S^4_q$ is not a quantum homogeneous space and so its differential geometry provides yet another example of the phenomenon where `quantisation removes degeneracy'.

Just as a Lie group imposes a differential structure upon each of its homogeneous spaces, so it is true of matrix quantum groups (under certain mild circumstances). In \S\ref{se:non-univ} we use Majid's theory of framed noncommutative manifolds \cite{Ma1,Ma:sphere} to equip the spaces $S^7_q$ and $\CP^3_q$ with finite-dimensional differential structures, finding some important differences from the classical case. It is a recognised fact in quantum group theory \cite{brmaj} that, if one wants to have a connection on a quantum principal bundle with non-universal calculi, then one needs to insist that the structure quantum group is equipped with a bicovariant differential calculus. In our case, doing so means that the quantum sphere $S^7_q$ carries an eight-dimensional differential calculus, whereas the geometry of twistor space turns out to be seven-dimensional, each being one higher than their classical value (a similar phenomenon occurs in \cite{bl:sph}). Nevertheless we find that the base space $S^4_q$ carries a well-defined, four-dimensional differential calculus that is obtained as the $\SU_q(2)$-invariant part of the differential calculus on the total space $S^7_q$. 

In \S\ref{se:bi} we introduce a Hodge duality structure for differential two-forms on the four-sphere $S^4_q$ and an associated set of anti-self-duality equations. Using the noncommutative differential geometry of the twistor fibration $\CP^3_q\to S^4_q$ we are able to construct an `instanton' solution to these equations. To our knowledge, this constitutes the first example of a genuine instanton connection on a $q$-deformed quantum four-sphere; it is hoped that this will lead to a more general description of instantons on noncommutative spheres in the context of $q$-deformed quantum groups. 

\emph{Notation}  
We refer to \cite{ks:qgr,Ma:book} for the axioms of unital Hopf $*$-algebras and basic notions such as modules and comodules. Given a Hopf algebra $H$ over the complex numbers $\C$, we denote its coproduct,
counit and antipode by $\Delta:H \rightarrow H\otimes H$,
$\ep:H\rightarrow \C$ and $S:H\to H$, respectively. We use Sweedler
notation for the coproduct, $\Delta h=h\o \otimes h\t$; iterated to $(\Delta \otimes \id) \circ \Delta h = (\id \otimes \Delta)
\circ \Delta h = h\o \otimes h\t \otimes h \thr$ and so on, with
summation inferred. If $\Delta_V:V\to V\otimes H$ is a right $H$-comodule, we also use the Sweedler-like notation $\Delta_V(v)=v\bz\otimes v\uo$ for each $v\in V$. When not at risk of confusion, we use the Einstein convention of summing over repeated indices. Finally, we use the notation $\sigma:=q-\iq$.

\subsubsection*{Acknowledgments}
Both authors were partially supported by the Italian project
`Cofin08--Noncommutative Geometry, Quantum Groups and Applications'. They thank an anonymous referee for pointing out some errors in an earlier version of the paper and are grateful to Francesco D'Andrea, Gaetano Fiore and Chiara Pagani for several suggestions and improvements. SJB acknowledges support from FNR Luxembourg under the PDR-project 894130 and from the EU project `Geometry and Symmetry of Quantum Spaces' PIRSES--GA-2008-230836.

\section{Noncommutative Differential Structures}\label{se:prelims}

In this paper we study differential structures on quantum principal bundles. Since the basic concepts 
are by now surely well-known, we 
confine ourselves to an easy review while occasionally referring to the known literature.

\subsection{Differential calculi on quantum groups}\label{se:dcqg}
Let $(P, m_P)$ be a unital $*$-algebra. Recall that a {\em first order differential calculus} over
$P$ is a pair $(\Omega^1P,\D)$, where $\Omega^1P$ is a
$P$-$P$-bimodule (the one-forms) and $\D:P
\rightarrow \Omega^1P$ is a linear map obeying the Leibniz rule and such that $\Omega^1P$ is spanned by elements of the form $a\, \D b$. A differential calculus $(\Omega^1P,\D)$ is said to be a {\em $*$-calculus} if it obeys $(\D p)^*=\D(p^*)$ for all $p\in P$.

The  {\em universal differential calculus} $(\widetilde\Omega^1P,\widetilde\D)$ over $P$ is the $P$-$P$-bimodule $\widetilde\Omega^1P:=\textrm{Ker}\,m_P$ equipped with the differential $\widetilde \D p=1\otimes p-p\otimes 1$ for each $p\in P$. It is universal in that any differential calculus $(\Omega^1P,\D)$ over $P$ arises as the quotient of $(\widetilde\Omega^1P,\widetilde\D)$ by some $P$-$P$-sub-bimodule $N_P \subset \textrm{Ker}\,m_P$.

Let $H= (H,\Delta,S,\ep_H)$ be a Hopf algebra. 
A calculus on $H$ is {\em left-covariant} if the coproduct $\Delta$,
as a left coaction of $H$ on itself, extends to a left
coaction on one-forms, $\Delta_L:\Omega^1H \rightarrow H\otimes\Omega^1H$ such
that $\D$ is an intertwiner and $\Delta_L$ is a bimodule map. In this situation, one knows that $\Omega^1 H\simeq H\otimes\Lambda^1$, where $\Lambda^1$ is the vector space of left invariant one-forms, and $\Lambda^1\simeq H^+/I_H$, for $I_H$ 
a right ideal of the augmentation ideal $H^+:=\textup{Ker}\,\ep_H$.

More explicitly, this classification is given by the linear isomorphism
\begin{equation}\label{iso-cl}
\varpi:H\otimes H^+\to \widetilde\Omega^1H,\qquad \varpi(h\otimes g):=hS(g\o)\otimes g\t.
\end{equation}
Given a right ideal $I_H$ of $H^+$, the vector space $H\otimes I_H$ is carried by $\varpi$ onto a left-covariant $H$-$H$-sub-bimodule $N_H$ of $\widetilde\Omega^1H$ and every such bimodule $N_H$ arises in this way \cite{wor}. This also shows that $(\Omega^1H,\D)$ is a $*$-calculus if and only if the corresponding ideal $I_H$ is such that $S((I_H)^*)\subseteq I_H$.

Right covariant differential calculi over $H$ are defined similarly. A calculus is called {\em bicovariant} if its both left and right covariant. Define the right adjoint coaction by
$$
\textup{Ad}_R:H\to H\otimes H,\qquad
\textup{Ad}_R(h)=h\t\otimes S(h\o) h\thr.
$$
Then, a left covariant calculus is bicovariant if and only if the corresponding ideal $I_H$ of $H^+$ is $\textup{Ad}_R$-stable, meaning $\textup{Ad}_R(I_H)\subseteq I_H\otimes H$. 

A collection of examples of bicovariant differential calculi on Hopf algebras, each playing an important role in the present paper, are given in \S\ref{se:su2}. 

\subsection{Quantum principal bundles}\label{qpb}\label{se:frame} 
Next, we recall the notion of a quantum principal bundle and its
associated differential structures. Let $(P,H,\delta_R)$ be a triple consisting of a unital $*$-algebra $P$, a unital Hopf $*$-algebra $H$ and a right coaction $\delta_R:P \rightarrow P \otimes H$ which makes $P$ into a right $H$-comodule $*$-algebra. We define the subalgebra 
$$
M:=P^H=\{p \in P: \delta_R(p)=p\otimes 1\},
$$
consisting of elements which are invariant under the coaction.
\begin{defn}
A {\em quantum principal bundle} is a triple $(P,H,\delta_R)$ as above such that the following canonical map is bijective:
\begin{equation}\label{can map} 
\chi:P\otimes_MP\to
P\otimes H, \qquad p\otimes_M p'\mapsto p\,\delta_R(p').
\end{equation}
\end{defn}

Bijectivity of the canonical map \eqref{can map} is sufficient to have a principal bundle at the `universal level', 
{\em i.e.} in the case where the algebras $P$, $H$ $M$ are all equipped with their universal differential calculi 
\cite{brmaj}. For a principal bundle with non-universal calculi, one needs to impose extra conditions which guarantee compatibility of the differential
structures on the total space $P$ and on the structure quantum group $H$. 

Indeed, suppose that $P$, $M$ are equipped with differential calculi $\Omega^1P$, $\Omega^1M$ defined by sub-bimodules $N_P$, $N_M$ respectively, and that $H$ carries a left covariant differential calculus $\Omega^1H$ defined by a right ideal $I_H$. Compatibility of these differential structures means that the calculi should satisfy the three conditions
\begin{equation}
N_M = N_P \cap \widetilde \Omega^1M, \qquad 
\delta_R(N_P)\subseteq N_P\otimes H,\qquad 
\textup{ver}(N_P)=P\otimes I_H, \label{submodules}
\end{equation}
where $\textup{ver}(p \otimes p')=p\,\delta_R(p')$ is the canonical map which generates the vertical one-forms. 

The role of the first condition is to ensure that $\Omega^1M$ is spanned by
elements of the form $m \, \D n$ with $m,n \in M$ and is hence obtained
by restricting the calculus on $P$. The second condition in \eqref{submodules}
is sufficient to ensure covariance of
$\Omega^1P$. Finally, the third condition ensures that the map 
$$
\textup{ver}:\Omega^1P \rightarrow P\otimes \Lambda^1, \qquad
\Lambda^1 \simeq H^+/I_H,
$$
is well-defined and yields exactness of the following sequence \cite{Ma:sphere}:
\begin{equation}\label{eqn exact seq non-univ}0 \rightarrow
P(\Omega^1M)P \rightarrow \Omega^1P \xrightarrow{\textup{ver}}P\otimes\Lambda^1 \rightarrow 0 .
\end{equation}
 
A special type of quantum principal bundle of later use in the paper comes from the following construction  \cite{brmaj}.
Suppose that $P$ is itself a Hopf algebra equipped with a surjection of Hopf
algebras $\pi:P \rightarrow H$. Then there is a right coaction of $H$
on $P$ by coproduct and projection to $H$,
$$
\delta_R:P \rightarrow P \otimes H,\qquad
\delta_R:=(\id\otimes\pi)\circ\Delta.
$$
The base algebra $M=P^H$ of coinvariants is now called a {\em quantum homogeneous space}. 
We do not dwell here upon the extra conditions needed for the canonical map as in \eqref{can map} be bijective 
so that $(P,H,\delta_R)$ is a quantum principal bundle with universal differential calculi. Rather we shall do this for the particular cases in which we are interested.

Next, suppose the Hopf algebras $P$ and $H$ carry left-covariant differential calculi $\Omega^1P$ and
$\Omega^1H$. They are respectively defined by
right ideals $I_P$ of $P^+$ and $I_H$ of $H^+$. For a quantum principal bundle with non-universal calculi we 
need the compatibility conditions
\begin{equation} \label{hom bundle}
(\textup{id}\otimes\pi)\textup{Ad}_R(I_P) \subset I_P\otimes H,
\qquad \pi(I_P)=I_H.
\end{equation}
A choice of left-covariant calculus on $P$ satisfying these
conditions automatically gives a principal bundle with non-universal
calculi \cite{Ma:sphere}.

Just as in the classical case, one introduces the notion of a vector bundle associated to a quantum principal bundle. For this, we need not only
a quantum principal bundle $\delta_R:P \rightarrow P \otimes H$ as
above but also a right $H$-comodule $\Delta_R:V\to V\otimes H$. 
The role of the space of sections of the associated vector bundle is played by the following object.

\begin{defn}
Let $(P,H,\delta_R)$ be a quantum principal bundle and let $V$ be a right $H$-comodule. Then the {\em associated vector bundle} $\mathcal{M}(V)$ is the vector space
\begin{equation}\label{assoc-bd}
\mathcal{M}(V):=(P \otimes V)^H
\end{equation}
of coaction invariant elements, where the
vector space $P\otimes V$ is equipped with the right tensor product
coaction. By construction, $\mathcal{M}(V)$ is an $M$-bimodule.
\end{defn}

To give the notion of a framing on a quantum space $M$, we also require a `soldering form' $\theta:V
\rightarrow P(\Omega^1M)$ for which the induced left $M$-module map
\begin{equation}\label{indiso}
s_\theta: \mathcal{M}(V) \rightarrow \Omega^1M, \quad p\otimes v \mapsto
p\,\theta (v)
\end{equation}
is an isomorphism. This leads to the following definition.

\begin{defn}
An algebra $M$ is said to be a {\em framed quantum manifold} if it
is the base of a quantum principal bundle,
$M=P^H$, to which the differential calculus $\Omega^1M$ is an associated vector bundle equipped with a soldering form.
\end{defn}

In the case where $M$ is a quantum homogeneous space, if the
conditions in eq.~\eqref{hom bundle} are satisfied then the algebra
$M=P^H$ is automatically framed by the bundle $(P,H,\delta_R)$
\cite{Ma:sphere}. The $H$-comodule $V$ and soldering form $\theta$
are given explicitly by the formul{\ae}
\begin{equation}\label{frame}
V=(P^+ \cap M)/(I_P \cap M), \qquad \Delta_Rv= \tilde v\t\otimes S\pi(\tilde v\o), \qquad\theta(v)=S(\tilde v\o)\D\tilde v\t
\end{equation}
for any representative $\tilde v$ of $v$ in $P^+ \cap M$, where
$\Delta(\tilde v)=\tilde v\o\otimes \tilde v\t$ is the
coproduct on $P$.

\subsection{The quantum group $\SU_q(2)$}\label{se:su2}
With $0<q<1$ a real deformation parameter, the algebra
$\A[\SU_q(2)]$ of coordinate functions on the quantum group
$\SU_q(2)$ is the associative unital algebra generated by the entries
of the matrix
\begin{equation}\label{m-gens}
\mathbf{a}=(\sa_i{}^j)=\begin{pmatrix}a&b\\c&d\end{pmatrix}
\end{equation}
obeying the relations
\begin{eqnarray}\label{mat rels}
& ab=qba,\qquad ac=qca,\qquad bd=qdb,\qquad cd=qdc,\\\nonumber &
bc=cb,\qquad ad-da=(q-\iq)bc .
\end{eqnarray}
There is also the determinant relation
$$
ad-qbc=1
$$
which, with the last equation in \eqref{mat rels}, implies in addition that
$$
da-\iq bc=1.
$$
The algebra $H:=\A[\SU_q(2)]$ has a matrix coproduct and counit, defined on
generators by $\Delta(\sa_i{}^j)=\sa_i{}^\mu\otimes\sa_\mu{}^j$ and
$\ep(\sa_i{}^j)=\delta_i^j$, and extended as algebra maps. 
Using the $R$-matrix
\begin{equation}\label{R mat}R=(R_i{}^j{}_k{}^l)=\begin{pmatrix}q&0&0&0\\0&1&0&0\\0&q-\iq&1&0\\0&0&0&q\end{pmatrix},
\end{equation}
the algebra relations \eqref{mat rels} may be written more concisely as
\begin{equation}\label{c-rels}
 R_i{}^\alpha{}_k{}^\beta \sa_\alpha{}^j
\sa_\beta{}^l = \sa_k{}^\beta\sa_i{}^\alpha R_\alpha{}^j{}_\beta{}^l,
\end{equation}
for $i,k,j,l=1,2$, or with the `compact matrix notation' of \cite{Ma:book},
in a shorthand expression
$$
R\, \bfa_1\bfa_2=\bfa_2\bfa_1 R .
$$
Moreover, $\A[\SU_q(2)]$ is equipped with an anti-linear involution
\begin{equation}\label{star}
\mathbf{a}^*=\begin{pmatrix}a^*&b^*\\c^*&d^*\end{pmatrix}:=\begin{pmatrix}d&-q
c\\-\iq b&a\end{pmatrix}
\end{equation}
and an antipode defined by $S(\sa_k{}^l)=(\sa_l{}^k)^*$, both extended as anti-algebra maps. These structures together make $H=\A[\SU_q(2)]$ into a Hopf $*$-algebra. The latter has a
coquasitriangular structure given on
generators by
\begin{equation}\label{ct}
\mathcal{R}:H\otimes H\to \C,\qquad \mathcal{R}(\sa_i{}^j\otimes \sa_k{}^l)=\zeta \, R_i{}^j{}_k{}^l
\end{equation}
and extended as a Hopf bicharacter. The factor
$\zeta=q^{-1/2}$ here is a `quantum group normalisation', which
ensures that $\mathcal{R}$ is compatible with the antipode in the
sense that
$$
\mathcal{R}(h\otimes Sg)=\mathcal{R}^{-1}(h\otimes
g)=\mathcal{R}(Sh\otimes g) , \qquad 
\textup{for all} \quad h,g\in \A[\SU_q(2)] .
$$ 
The two-dimensional vector space
$V_{\frac{1}{2}}:=\textrm{Span}_\C\{x,y\}$ is the {\em fundamental comodule} for $\A[\SU_q(2)]$
and is equipped with the right coaction
\begin{equation}\label{fun-co}
\Delta_{\frac{1}{2}}:V_{\frac{1}{2}}\to V_{\frac{1}{2}}\otimes \A[\SU_q(2)],\qquad
\begin{pmatrix}x&y\end{pmatrix}\mapsto \begin{pmatrix}x&y\end{pmatrix}\otimes
\begin{pmatrix}a&-qc^*\\c&a^*\end{pmatrix}.
\end{equation}
Let $\A[\C^2_q]$ denote the unital algebra generated by $x,y$
subject to the relation $xy=qyx$.
Then by extending the map \eqref{fun-co} as an
algebra map, the algebra $\A[\C^2_q]$ becomes a right
$\A[\SU_q(2)]$-comodule algebra,
\begin{equation}\label{ext-co}
\A[\C^2_q]\to\A[\C^2_q]\otimes\A[\SU_q(2)].
\end{equation}
For each $j=0, \tfrac{1}{2},1, \tfrac{3}{2},\ldots$, let $V_j$ denote the
$2j+1$-dimensional vector space spanned by the set of
polynomials in $\A[\C^2_q]$ of degree $2j$. Then by restricting the
coaction \eqref{ext-co}, we obtain a right comodule structure
$$
\Delta_j:V_j\to  V_j\otimes\A[\SU_q(2)].
$$
The index $j$ is called the {\em spin} of the corepresentation. Each
of these comodules $V_j$ is known to be an irreducible unitary corepresentation
of the Hopf algebra $\A[\SU_q(2)]$ and every such finite-dimensional corepresentation
arises in this way \cite{ks:qgr}.

We shall also need the `co-opposite' quantum group $\SU^\cop_q(2)$, whose $*$-algebra of coordinate functions
$\tH:=\A[\SU^\cop_q(2)]$  is just a copy of $H:=\A[\SU_q(2)]$. As for the coalgebra structure 
$(\widetilde\ep, \tS, \tDelta)$, one keeps $\widetilde\ep=\ep$ while $\tS = S^{-1}$ and the coproduct is changed to  
$$
\tDelta(a_i{}^j)= a_\mu{}^j\otimes a_i{}^\mu. 
$$
The two-dimensional vector space
$\tV_{\frac{1}{2}}:=\textrm{Span}_\C\{\tilde x,\tilde y\}$, the fundamental comodule for $\A[\SU^\cop_q(2)]$, 
now carries the right coaction
\begin{equation}\label{fun-co2}
\Delta_{\frac{1}{2}}:\tV_{\frac{1}{2}}\to \tV_{\frac{1}{2}}\otimes \A[\SU^\cop_q(2)], \qquad
\begin{pmatrix}\tilde x&\tilde y\end{pmatrix}\mapsto \begin{pmatrix}\tilde x&\tilde y\end{pmatrix}\otimes
\tS \begin{pmatrix} a& -q c^*\\ c& a^*\end{pmatrix}.
\end{equation}
As before, for each $j=0, \tfrac{1}{2},1, \tfrac{3}{2},\ldots$, there will be corepresentations 
$$
{\widetilde{\Delta}}_j:\tV_j\to \tV_j\otimes \A[\SU^\cop_q(2)].
$$ 
the spin $j$ irreducible comodule $\tV_j$ being the $2j+1$-dimensional vector space spanned by the collection 
of polynomials
of degree $2j$ within the unital algebra $\widetilde{\A}[\C^2_q]$ generated 
by the elements $\tilde x,\tilde y$,
subject to the relation $\tilde x \tilde y = q \tilde y\tilde x$.

Finally, we shall need the classical subgroup $\U(1)$ of the quantum
group $\SU_q(2)$. The algebra $H'=\A[\U(1)]$ of coordinate functions
on the group $\U(1)$ is the commutative unital $*$-algebra generated
by the mutually conjugate elements $t$ and $t^*$, subject to the relations $tt^*=t^*t=1$. It becomes a
Hopf algebra when equipped with the coproduct, counit and antipode
defined on generators by
$$
\Delta(t)=t\otimes t,\qquad \ep(t)=1,\qquad S(t)=t^*;
$$
as usual, the coproduct and counit are extended as $*$-algebra maps, the
antipode as a $*$-anti-algebra map. There is a canonical Hopf
$*$-algebra surjection
\begin{equation}\label{hopfproj}
\pi:\A[\SU_q(2)]\to\A[\U(1)], \qquad
\pi\begin{pmatrix}a&-qc^*\\c&a^*\end{pmatrix}=\begin{pmatrix}t&0\\0&t^*\end{pmatrix},
\end{equation}
which, in the classical case, is the dual of the group inclusion $\U(1)\hookrightarrow\SU(2)$.

The irreducible $\A[\U(1)]$-comodules are all one-dimensional and labeled by an
integer $k\in\ZZ$. We write $U_k$ for the irreducible comodule
spanned by the vector $u_k$ with coaction 
$$
\Delta_k':U_k\to  U_k\otimes H', \qquad \Delta_k'(u_k)=
u_k\otimes t^k,
$$ 
where $t$ is the generator of the algebra $\A[\U(1)]$ and we define $t^{-k}:=(t^*)^k$ for each $k\in\mathbb{N}$.

\begin{example}\label{ex:4d}
It is known \cite{wor} that there is no three-dimensional bicovariant differential calculus on $\SU_q(2)$, whence the need for a four-dimensional calculus for bicovariance. This fact will play a very important role in the geometry we discuss in the present paper.

The 4D${}_+$ differential
calculus $\Omega^1\SU_q(2)$ on $H=\SU_q(2)$, originally described
in \cite{wor}, is defined by the right ideal $I_H$
of $\textrm{Ker}\,\ep_H$ generated by the nine elements
\begin{equation}\label{nine}
b^2,~c{}^2,~b(a-d),~c(a-d),~a^2+q^{2}d^2-(1+q^2)(ad+\iq bc),
\end{equation}
$$\st\,b,~\st\,c,~\st\,(a-d),~\st\,(q^2a+d-(q^2+1)),$$
where $\st:=q^2a+d-(q^3+q^{-1})$. From the discussion above, we know that
this ideal determines a left-covariant $*$-calculus on the Hopf algebra
$H=\A[\SU_q(2)]$. The ideal $I_H$ is stable under the adjoint coaction $\textrm{Ad}_R$, so that the 
calculus is bicovariant. 
\end{example}

It is straightforward to check that the vector space
$\Lambda^1\simeq H^+/I_H$ of left-invariant
one-forms in the calculus $\Omega^1\SU_q(2)$ is four-dimensional and spanned by the elements $b$, $c$
and the elements $a_0$, $a_z$ defined by the equations
\begin{align}\label{4dspan}
a-1=(\iq-1)a_0+(q-1)a_z, \qquad
d-1=(q-1)a_0+(\iq-1+\sigma^2\iq)a_z.
\end{align}

\begin{example}\label{ex:4dopp}
From this example, one finds that the co-opposite algebra $\tH=\A[\SU^\cop_q(2)]$ is also equipped with a four-dimensional first order differential calculus.
Indeed, since the $*$-algebras $\A[\SU_q(2)]$ and $\A[\SU^\cop_q(2)]$ are isomorphic, we find a left-covariant calculus $\Omega^1\SU_q^\cop(2)$ determined by the right ideal $I_{\tH}$ of $\tH$ which is a copy of the ideal $I_{H}$ generated as in 
\eqref{nine}. 
With the `new' coproduct $\tDelta$ on $\tH$, the notions of left and right covariance are of course co-opposite to those of $H$, although it is nevertheless straightforward to deduce that the ideal $I_{\tH}$ is stable under the right adjoint coaction $\textrm{Ad}_R$ of $\tH$, so that the calculus is bicovariant. As bimodules, the calculi $\Omega^1\SU_q(2)$ and $\Omega^1\SU_q^\cop(2)$ are the same.
\end{example}

Our final example concerns a quantum differential calculus on the classical group $\U(1)$ described as before 
by the commutative Hopf algebra $H'=\A[\U(1)]$.
\begin{example}
From the calculus $\Omega^1\SU_q(2)$ we immediately obtain a differential structure $\Omega^1\U(1)$ on the
classical Hopf algebra $H'=\A[\U(1)]$ 
in terms of an $\textrm{Ad}_R$-stable right ideal
$I_{H'}:=\pi(I_H)$ of $H'{}^+$, where $\pi:H\to H'$ is the projection \eqref{hopfproj}. The ideal $I_{H'}$ is generated
by the three elements
\begin{equation}\label{onegens}
t{}^2+q^{2}t^*{}^2-(1+q^2),\quad \st\,(t-t^*),\quad(q^2t+t^*-(q^2+1))\, \st,
\end{equation}
again with $\st=q^2t+t^*-(q^3+\iq)$, where $t,t^*$ are the
generators of $H'$.

In this case, one easily finds (as in \cite{bl:sph} for example) that the vector space
$(H')^+/I_{H'}$ of left-invariant one-forms is one-dimensional and spanned by the element $t-1$. Since the defining ideal is again $\textup{Ad}_R$-stable, the calculus on $\U(1)$ is also bicovariant.
\end{example}

\section{The Quantum Hopf Fibration}\label{se:hopf-fib}
The classical $\SU(2)$ Hopf fibration is nothing other than the
canonical quaternionic line bundle over the projective space
$\HH\mathbb{P}^1$ ({\em cf}. \cite{ma:gymf} for information relevant to the present
paper). In this section, we keep this quaternionic interpretation
firmly in mind and use it to construct a deformed version of the
Hopf fibration.

\subsection{Deformations of quaternionic spaces}
We shall
need then a $q$-version of the vector space $\HH$ of quaternions. We
begin with a $q$-deformed analogue of the Euclidean space $\RR^4$,
which we later equip with a quaternionic structure.
The algebra $\A[\RR^4_q]$ we use was proposed in \cite{ma:twist} as giving
a natural $q$-analogue of the Euclidean space $\RR^4$, since it
possesses a natural `metric' with the correct Euclidean signature.
It coincides with the one used in \cite{fiore,fiore2}.

\begin{defn}
The coordinate algebra $\A[\RR^4_q]$ is the unital $*$-algebra
generated by the entries of the matrix
\begin{equation}\label{quatgens}
\bfx:=(\sx_i{}^j)=\begin{pmatrix}\iq z_1&-z_2^*\\z_2&z_1^*\end{pmatrix}
\end{equation}
subject to the relations
\begin{eqnarray}\label{quat rels}
& z_1z_2=\iq z_2z_1 \quad \textup{and} \quad z_2^*z_1^*=\iq z_1^*z_2^*, 
\quad z_1z_2^*=qz_2^*z_1 \quad \textup{and} \quad z_2z_1^*=qz_1^*z_2, \nonumber  \\
& z_1z_1^*=z_1^*z_1 \quad \textup{and} \quad z_2^*z_2-z_2z_2^*= (1-q^{-2})z_1z_1^*.
\end{eqnarray}
\end{defn}
\noindent
Notice the similarities and differences with the relations \eqref{mat rels}.
They can be written more compactly, for $i,k,j,l=1,2$, as
\begin{equation}\label{quat rels2}
  R_k{}^\beta{}_i{}^\alpha \sx_\alpha{}^j
\sx_\beta{}^l =  \sx_k{}^\beta\sx_i{}^\alpha R_\alpha{}^j{}_\beta{}^l,
\qquad \textup{or} \qquad R_{21}\, \bfx_1\bfx_2=\bfx_2\bfx_1R \, ,
\end{equation}
where the $R$-matrix $R_{21}$ is defined in terms of the one in \eqref{R mat} by 
\begin{equation}\label{R2}
(R_{21})_i{}^j{}_k{}^l=R_k{}^l{}_i{}^j, \qquad i,k,j,l=1,2.
\end{equation}

In analogy with the approach of \cite{lprs:ncfi}, we introduce a quaternionic involution $\mathrm{J}_q$ on $\A[\RR^4_q]$, defined on
the generators \eqref{quatgens} by
\begin{equation}\label{J-quat}
\mathrm{J}_q:\A[\RR^4_q]\to\A[\RR^4_q],\qquad
\mathrm{J}_q\begin{pmatrix}\iq z_1&-z_2^*\\z_2&z_1^*\end{pmatrix}=q^{-1/2}\begin{pmatrix}z_2^*&z_1\\-
z_1^*&q z_2\end{pmatrix}
\end{equation}
and extended as an anti-algebra map. One readily checks
that $\mathrm{J}^2_q=-\id$, as claimed.
The map \eqref{J-quat} equips $\A[\RR^4_q]$ with a quantum analogue of a `quaternionic
structure' (on the noncommutative space $\RR^4_q$). Indeed, in the classical limit, we would be identifying
the set of quaternions $\HH$ with the set of complex $2\times 2$
matrices of the form
$$
c_1+c_2j\in\HH \mapsto \begin{pmatrix}c_1&-\bar c_2\\ c_2& \bar
c_1\end{pmatrix},
$$
with the map $\mathrm{J}_q$ corresponding in the limit to
right multiplication by the quaternion $j$.

\begin{defn}\label{def:quats}
We define the coordinate algebra $\A[\HH_q]$ of the $q$-deformed
quaternions $\HH_q$ to be the $*$-algebra $\A[\RR^4_q]$ having in addition 
the quaternionic involution $\mathrm{J}_q$ of eq.~\eqref{J-quat}.
\end{defn}

The algebra $\A[\HH_q]$ carries canonical commuting
coactions defined by
\begin{equation}\label{lc}
\delta_R :\A[\HH_q]\to\A[\HH_q]\otimes\A[\SU_q(2)],  \qquad
\sx_i{}^j \mapsto \sx_i{}^\mu\otimes \sa_\mu{}^j,
\end{equation}
for $i,j=1,2$, as well as 
\begin{equation}\label{rc}
\delta_L :\A[\HH_q]\to\A[\SU^\cop_q(2)]\otimes\A[\HH_q], \qquad
\sx_i{}^j \mapsto \tS(a_i{}^\mu) \otimes\sx_\mu{}^j
\end{equation}
for $i,j=1,2$, and both extended as $*$-algebra maps, where $\A[\SU^\cop_q(2)]$ denotes the co-opposite Hopf algebra of $\A[\SU_q(2)]$.

Using these, 
from the `one-dimensional'
quaternionic space we pass to its two-dimensional analogue $\A[\HH^2_q]:=\A[\HH_q]\otimes\A[\HH_q']$, defined as the 
(braided) tensor product of two copies of the algebra $\A[\HH_q]$. The cross-relations in the tensor product are obtained by requiring them to be covariant under the right coaction \eqref{lc} {\em i.e.} by constructing a
braided tensor product algebra in the category of right $\A[\SU_q(2)]$-comodules  \cite{ma:twist}.

\begin{defn}We define $\A[\HH^2_q]$ to be the (braided) tensor product $\A[\HH_q]\otimes\A[\HH_q']$ generated by two copies of
$\A[\HH_q]$, whose generators we denote by
$$
\mathbf{x}=\begin{pmatrix}\iq z_1&-z_2^*\\z_2&z_1^*\end{pmatrix},
\qquad \mathbf{x}'=\begin{pmatrix}\iq z_3&-z_4^*\\
z_4&z_3^*\end{pmatrix}
$$
for $\A[\HH_q]$ and $\A[\HH_q']$ respectively, with commutation relations as above in \eqref{quat rels}. 
With $\zeta=q^{-1/2}$ the normalisation as in \eqref{ct}, the cross-relations between $\mathbf{x}$ and $\mathbf{x}'$ in the algebras $\A[\HH_q]$ and $\A[\HH_q']$ are found to be
\begin{align}\label{quat rels3}
z_1z_3   & =\zeta q \, z_3z_1, \qquad  z_1^*z_3   =\zeta \, z_3z_1^*,   \qquad
z_2z_3   = \zeta q \, z_3z_2,  \qquad   z_2^*z_3   =\zeta \, z_3z_2^*,  \nn \\ 
z_4z_1   & = \zeta q \, z_1z_4, \qquad  z_1z_4^* -(q-\iq)\, z_2^*z_3  =\zeta \, z_4^*z_1, \nn \\
z_4z_2   & =\zeta q \, z_2z_4,   \qquad  z_2z_4^*+ (1-q^{-2})\, z_1^*z_3  =\zeta \, z_4^*z_2,
\end{align}
together with the conjugate relations obtained using the $*$-structure.
\end{defn}

For later use, we observe that the cross-relations \eqref{quat rels3} may also be written as
\begin{align}
\label{quat rels4} \sx_1{}^i\sx'_1{}^k=\zeta\sx'_1{}^\beta\sx_1{}^\alpha R_\alpha{}^i{}_\beta{}^k,& \qquad \sx'_2{}^i\sx_1{}^k=\zeta\sx_1{}^\beta\sx'_2{}^\alpha R_\alpha{}^i{}_\beta{}^k, \\
\sx_2{}^i\sx'_1{}^k=\zeta\sx'_1{}^\beta\sx_2{}^\alpha R_\alpha{}^i{}_\beta{}^k,& \qquad \sx'_2{}^i\sx_2{}^k=\zeta\sx_2{}^\beta\sx'_2{}^\alpha R_\alpha{}^i{}_\beta{}^k, \nn
\end{align}
in terms of the $R$-matrix \eqref{R mat}. 

\begin{lem}
The algebra $\A[\HH^2_q]$ is made into a right $\A[\SU_q(2)]$-comodule algebra by a right coaction $\delta_R:\A[\HH^2_q]\to \A[\HH^2_q]\otimes \A[\SU_q(2)]$ defined on generators by
\begin{equation}\label{su coact}
\delta_R\left( \sx_i{}^j\otimes \sx'_k{}^l \right) = 
\left(\sx_i{}^\mu\otimes\sx'_k{}^\nu\right)\otimes \sa_\mu{}^j  \sa_\nu{}^l \, ,
\end{equation}
and extended as a $*$-algebra map.
\end{lem}

\proof We already know  that this coaction makes
$\A[\HH_q]$ and $\A[\HH_q']$ into right $\A[\SU_q(2)]$-comodule
$*$-algebras. It is straightforward to check that the cross-relations in the tensor product
$\A[\HH^2_q]$ are covariant as well for it
(they were indeed derived having this property in mind). 
As said, all of the above amounts to saying that the
algebra $\A[\HH_q]\otimes\A[\HH_q']$ is the braided tensor product
algebra in the category of right $\A[\SU_q(2)]$-comodules.
\endproof

Finally, the quaternionic involution $\mathrm{J}_q$ extends to $\A[\HH^2_q]$
naturally according to 
\begin{equation}\label{eqn J}
\mathrm{J}_q:\A[\HH^2_q]\to\A[\HH^2_q], \qquad
\mathrm{J}_q\begin{pmatrix}\iq z_1&-z_2^*\\z_2
& z_1^*\\\iq z_3&-z_4^*\\z_4&z_3^*\end{pmatrix}=q^{-1/2}\begin{pmatrix}z_2^*&
z_1\\-z_1^* &qz_2\\z_4^*&z_3\\-z_3^*&qz_4\end{pmatrix},
\end{equation}
extended as a $*$-anti-algebra map. Again one finds that
$\mathrm{J}_q^2=-\id$, as it should be.

\subsection{The quantum Hopf bundle}\label{hopf fib}
We are ready to  construct a quantum seven-sphere $\A[S^7_q]$ and
equip it with a right coaction of the quantum group $\A[\SU_q(2)]$,
thus yielding a quantum principal bundle whose `base space'  is a quantum four-sphere $\A[S^4_q]$.

\begin{lem}
Define $\det_q \mathbf{x}=q^{-2}z_1z_1^*+z_2^*z_2$ and $\det_q
\mathbf{x}'=q^{-2}z_3z_3^*+z_4^*z_4$. Then
$$r^2:=\det\nolimits_q\mathbf{x}+\det\nolimits_q\mathbf{x}'$$
is a central element of the algebra $\A[\HH^2_q]$.\end{lem}

\proof The element $\det_q \mathbf{x}$ is known to be central in the
subalgebra generated by $\mathbf{x}$, and it is easy to check that
$\det_q \mathbf{x}$ commutes with the generators $\mathbf{x}'$. Similarly, 
$\det_q \mathbf{x}'$ is central in the subalgebra generated by $\mathbf{x'}$ 
and commutes with the generators $\mathbf{x}$.
\endproof

\begin{defn} The coordinate algebra $\A[S^7_q]$ of the
quantum sphere $S^7_q$ is the quotient of the algebra
$\A[\HH^2_q]$ by the two-sided ideal generated by the central
element $r^2-1$.\end{defn}

The ideal generated by $r^2-1$ is preserved by the right
coaction \eqref{su coact}, whence the latter descends to a coaction of
$\A[\SU_q(2)]$ on the algebra $\A[S^7_q]$, given by the
same formula.

\begin{rem}\label{rem:cp}\textup{The sphere $S^7_q$ has many classical points, 
amongst which the most obvious one is the one corresponding to the character 
$\phi:\A[S^7_q]\to\C$ which maps $z_4\mapsto 1$, $z_4^*\mapsto 1$ and the 
other generators to zero. We shall use this particular character later on ({\em cf}. Prop.~\ref{isom1}), when we
come to consider the quantum homogeneous space structure of
$S^7_q$.}\end{rem}

Next we come to the quantum Hopf fibration itself. On the free right
$\A[S^7_q]$-module $\mathcal{E}:=\C^4\otimes \A[S^7_q]$, there is the
canonical Hermitian structure
$h:\mathcal{E}\times\mathcal{E}\to\A[S^7_q]$ given by
$$
h(|\xi\ra,|\eta\ra):=\sum_{j=1}^4 (\xi_j)^*\eta_j,
$$
where $|\xi\ra, |\eta\ra\in \mathcal{E}$. Then to each element
$|\xi\ra\in \mathcal{E}$ there is an associated element $\la\xi|$ in
the dual module $\mathcal{E}^*$, defined by the non-degenerate
pairing
\begin{equation}\label{herm}
\la\xi|\eta\ra:=h(|\xi\ra,|\eta\ra).
\end{equation}
Using this construction, we note that the two columns $|\phi_1\ra$,
$|\phi_2\ra$ of the matrix
\begin{equation}\label{u matrix}
\su:=\begin{pmatrix}|\phi_1\ra &|\phi_2\ra
\end{pmatrix}=\begin{pmatrix}\iq z_1&-z_2^*\\z_2 &z_1^*\\\iq z_3&-z_4^*\\z_4&z_3^*\end{pmatrix}
\end{equation}
are orthonormal with respect to the Hermitian pairing $\la~|~\ra$,
in the sense that $\la\phi_i|\phi_j\ra=\delta_{ij}$ for $i,j=1,2$.
It follows that $\su^*\su=\mathbbm{1}_2$ and hence that the matrix
\begin{equation}\label{pr-dec}
\pp:=\su\su^*=|\phi_1\ra\la\phi_1|+|\phi_2\ra\la\phi_2|
\end{equation}
is a self-adjoint idempotent ({\em i.e.} a projection) in the matrix
algebra $\M_4(\A[S^7_q])$.

\begin{prop}\label{pr:proj}
The entries of the projection $\pp=\su\su^*$ generate a
subalgebra of $\A[S^7_q]$ which is a deformation of the  
algebra of coordinate functions on the four-sphere $S^4$.\end{prop}

\proof We explicitly compute the elements of the projection $\pp$
and their commutation relations. The diagonal elements are
\begin{align*}
\pp_{11}=q^{-2}z_1z_1^*+z_2^*z_2,\quad \pp_{22}=z_2z_2^*+z_1^*z_1 = \pp_{11} , \\
\pp_{33}=q^{-2}z_3z_3^*+z_4^*z_4,\quad \pp_{44}=z_4z_4^*+z_3^*z_3 = \pp_{33} ,
\end{align*}
the last equality in both lines coming from the commutation relations in $\A[S^7_q]$. 
Thanks to the sphere relation $r^2=1$ in $\A[S^7_q]$, together they satisfy the relation
$$
\pp_{11}+\pp_{22}+\pp_{33}+\pp_{44}=2,
$$
Thus, only one of the $\pp_{ii}$'s is independent and we write them in terms of $x_0:=(2\pp_{11}-1)$:
\begin{equation}\label{hopf inc1}
\pp_{11}=\tfrac{1}{2}(1+x_0)=\pp_{22},\qquad
\pp_{33}=\tfrac{1}{2}(1-x_0)=\pp_{44}.
\end{equation}
As in the classical case, the elements $\pp_{12}$, $\pp_{34}$
vanish,
$$\pp_{12}=\iq z_1z_2^*-z_2^*z_1=0,\qquad \pp_{34}=\iq z_3z_4^*-z_4^*z_3=0,$$ and the
remaining elements are given by
\begin{align}\label{hopf inc2}
\pp_{13}&=q^{-2}z_1z_3^*+z_2^*z_4,\qquad \pp_{14}=\iq z_1z_4^*-z_2^*z_3, \nonumber \\
\pp_{23}&=\iq z_2z_3^*-z_1^*z_4,\qquad \pp_{24}=z_2z_4^*+z_1^*z_3,
\end{align}
with $\pp_{ji}=\pp_{ij}^*$ when $j>i$.
Again using the commutation relations in $\A[S^7_q]$, it is straightforward to
check that only two of these are independent. We take the independent ones to be $\pp_{13}$ and $\pp_{14}$, finding
that $\pp_{23}=-\zeta^{-1}q\, \pp_{14}^*$ and $\pp_{24}=\zeta\, \pp_{13}^*$. We write
$x_1:=2\pp_{13}$ and $x_2:=2\pp_{14}$, so that the projection $\pp$
has the form
\begin{equation}\label{defproj}
\pp=\tfrac{1}{2}\begin{pmatrix}1+x_0&0&x_1&x_2\\0&1+x_0&-\zeta^{-1}q x_2^*&\zeta x_1^*\\x_1^*&-\zeta^{-1}qx_2&1-x_0&0\\x_2^*&\zeta x_1&0&1-x_0\end{pmatrix}.
\end{equation}
By construction $\pp^*=\pp$, so that $x_0^*=x_0$ and $x_1^*$,
$x_2^*$ are conjugate to $x_1$, $x_2$, hence the
notation. The fact that $\pp^2=\pp$ is the easiest way to compute
the relations between the generators: doing so yields $x_0$ to be central and
\begin{eqnarray}\label{sphere rels1}
& x_2x_1=q^2x_1x_2, \quad \textup{and} \quad x_1^*x_2^*=q^2 x_2^*x_1^*, 
\qquad x_2^*x_1=qx_1x_2^* , \quad \textup{and} \quad x_1^*x_2=q x_2x_1^* \nonumber \\
& x_1x_1^*+x_2x_2^*+x_0^2=1,\qquad
q^{-1}x_1^*x_1+q^3x_2^*x_2+x_0^2=1 .
\end{eqnarray}  
Of
course, these relations can also be computed directly from the
relations in $\A[S^7_q]$. When $q\to 1$, the algebra generated by the
entries $x_1,x_1^*,x_2,x_2^*,x_0$ reduces to the algebra of
coordinate functions on the classical four-sphere $S^4$, in which
case $\pp$ is a function on $S^4$ taking values in the collection of rank
two projections in $\M_4(\C)$.\endproof

\begin{defn}The coordinate algebra $\A[S^4_q]$ is the $*$-algebra
generated by the projection elements $x_0, x_1, x_1^*, x_2, x_2^*$ subject
to the commutation relations \eqref{sphere rels1}. \end{defn}

As already mentioned, the algebra $\A[S^7_q]$ is
covariant under the right coaction \eqref{su coact} of the quantum
group $\A[\SU_q(2)]$. This coaction relates the algebra $\A[S^7_q]$
to the algebra $\A[S^4_q]$, as the following proposition shows.

\begin{prop}The algebra $\A[S^4_q]$ is the algebra of coinvariants 
under the coaction $\delta_R:\A[S^7_q]\to \A[S^7_q]\otimes \A[\SU_q(2)]$ defined in eq.~\eqref{su coact}.
\end{prop}

\proof We need to show that $\A[S^4_q]=\{x\in\A[S^7_q]~|~\delta_R(x)=x\otimes 1\}$, 
with $\delta_R$ the coaction of \eqref{su coact}. It is easy to check that the generators of $\A[S^4_q]$ are
coinvariants. For example,
\begin{align*}
\delta_R(x_1)&=2\big(q^{-2}\delta_R(z_1)\delta_R(z_3^*)+\delta_R(z_2^*)\delta_R(z_4)\big)
\\&=2\big(q^{-2}z_1z_3^*\otimes(aa^*+q^2cc^*)+z_2^*z_4\otimes(c^*c+a^*a)\big)
\\&=2(q^{-2}z_1z_3^*+z_2^*z_4)\otimes 1=x_1\otimes 1,
\end{align*}
with the same result on other generators computed similarly. 
This shows that the whole algebra $\A[S^4_q]$ consists of
coinvariants. However, we also need to check that there are no coinvariants in $\A[S^7_q]$ 
other than elements of $\A[S^4_q]$. This follows from the reasoning used to prove a
similar result in \cite{lpr:qsb}. It is clear that elements
$w_1\in\{\iq z_1,z_2,\iq z_3,z_4\}$, respectively
$w_{-1}\in\{z_1^*,-z_2^*,z_3^*,-z_4^*\}$, are weight vectors of
weight 1, respectively -1, in the fundamental corepresentation of
$\A[\SU_q(2)]$. As a consequence, all coinvariants are of the
form $(w_1w_{-1}-qw_{-1}w_1)^n$ and when $n=1$ these are just the
generators of $\A[S^4_q]$.\endproof

Thus we have a canonical inclusion of algebras
$\A[S^4_q]\hookrightarrow\A[S^7_q]$. Using methods analogous to those of \cite{lpr:qsb}, 
this algebra extension is shown to be a quantum principal
bundle with structure quantum group $\SU_q(2)$, {\em i.e.} the associated canonical map \eqref{can map} is bijective.

\subsection{Noncommutative twistor space}
In classical geometry, the twistor space $\CP^3$ is obtained as a real
six-dimensional manifold by making the quotient of $S^7$ by a
certain action of $\U(1)$. We give a quantum version
of this {\em via} a coaction of the Hopf algebra $\A[\U(1)]$ on
the quantum sphere algebra $\A[S^7_q]$ and seeking the subalgebra of
coinvariants.

From the canonical projection \eqref{hopfproj} we immediately obtain
a right coaction $\delta_R'$ of $\A[\U(1)]$ on $\A[S^7_q]$, by
applying the coaction $\delta_R$ of \eqref{su coact} then projecting
to $\A[\U(1)]$:
\begin{equation}\label{twistor coact}
\delta_R':\A[S^7_q]\to\A[S^7_q]\otimes \A[\U(1)], \qquad
\delta_R':=(\id\otimes \pi)\delta_R ,
\end{equation}
where $\pi$ is the surjection in \eqref{hopfproj}. Equivalently, one imposes a $\ZZ$-grading on the algebra $\A[S^7_q]$ for
which the generators have degrees
\begin{equation}\label{grad}
\textup{deg}(z_j)=1,\qquad \textup{deg}(z_j^*)=-1, \qquad \textup{for} \qquad  j=1,\ldots,4 .
\end{equation}

\begin{defn}The coordinate function algebra $\A[\CP^3_q]$ of quantum twistor space $\CP^3_q$ 
is the subalgebra of $\A[S^7_q]$ made of coinvariants for the coaction \eqref{twistor coact}; equivalently
the subalgebra of overall degree zero with respect to the
$\ZZ$-grading \eqref{grad}.\end{defn}

One  checks that the algebra extension $\A[\CP^3_q]\hookrightarrow
\A[S^7_q]$ is a quantum principal bundle, meaning that the corresponding canonical 
map \eqref{can map} is bijective.

It is clear that elements $z_jz^*_l$, for $j,l=1,\dots 4$, generate the (whole and only) algebra $\A[\CP^3_q]$ of degree zero elements.
In parallel with what we did for the four-sphere in \eqref{defproj}, we assemble these generators into a projection. 
In the notation of eq.~\eqref{u matrix}, let us define
$$
\sv=|\phi_1\ra=\begin{pmatrix}\iq z_1&z_2&\iq z_3&z_4\end{pmatrix}^{\textup{tr}}.
$$
where ${}^\textup{tr}$ denotes matrix transposition. We already know that $\sv^*\sv=1$, whence $\qp:=\sv\sv^*$  is a
projection in the matrix algebra $\M_4(\A[S^7_q])$. Explicitly, it
works out to be
\begin{equation}\label{twistor
proj}\qp=|\phi_1\ra\la\phi_1|=\begin{pmatrix}q^{-2}z_1z_1^*&\iq z_1z_2^*&q^{-2}z_1z_3^*&\iq z_1z_4^*\\\iq z_2z_1^*&z_2z_2^*&\iq z_2z_3^*&z_2z_4^*
\\q^{-2}z_3z_1^*&\iq z_3z_2^*&q^{-2}z_3z_3^*&\iq z_3z_4^*\\\iq z_4z_1^*&z_4z_2^*&\iq z_4z_3^*&z_4z_4^*\end{pmatrix}. \end{equation}

\noindent
As mentioned before, as generators of the algebra $\A[\CP^3_q]$ we take the entries of the matrix $\qp=(\qp_{jl})$. The relations in $\A[\CP^3_q]$ are inherited from those of $\A[S^7_q]$ although, not needing them, we refrain from writing them out explicitly. 
The $*$-structure on $\A[\CP^3_q]$ is also inherited from that of $\A[S^7_q]$, {\em i.e.}
$ \qp_{jl}^*=\qp_{lj}$, for $j,l=1,\ldots,4$.

In the classical limit
$q\to 1$,  one recovers the fact that $\qp$ is the tautological rank
one projector-valued function on $\C^4$, which we think of as the
defining projector of $\CP^3$.

With $\mathrm{J_q}:\A[S^7_q]\to\A[S^7_q]$ the quaternionic map in
eq.~\eqref{eqn J}, one has 
$$
\mathrm{J}_q|\phi_1\ra =
-q^{-1/2}|\phi_2\ra, 
$$
so that using eq.~\eqref{u matrix} we obtain, at the level of generators,
the matrix sum
\begin{equation}\label{matsum}
\pp=|\phi_1\ra\la\phi_1|+|\phi_2\ra\la\phi_2|=\qp+q\mathrm{J}_q(\qp) .
\end{equation}
In parallel with \cite{bm:qtt} for the $\theta$-deformed case, this sets up an obvious algebra
inclusion 
\begin{equation}\label{tw-bund}
\eta:\A[S^4_q] \hookrightarrow \A[\CP^3_q],
\end{equation} 
which is a noncommutative analogue of the classical Penrose twistor fibration
$\CP^3\rightarrow S^4$, thus justifying our thinking of $\CP^3_q$ as the twistor space of the quantum four-sphere $S^4_q$.

\section{Quantum Symmetries of the Hopf Fibration}\label{se:symm}

The central part of the paper is devoted to finding a non-universal
differential structure on the quantum sphere
$\A[S^4_q]$. Contrary to
the classical situation, there is no canonical way to go about doing
this. In order to simplify our task, we look for a differential structure which is covariant under a quantum group of symmetries of the Hopf fibration. In this section we explicitly construct such a quantum group.

\subsection{Quantum symmetries of $\HH^2_q$}To obtain a quantum
group of symmetries of the Hopf fibration, we begin by asking how
the quantum space $\HH^2_q$ behaves under linear transformations.
Having already used the right $\A[\SU_q(2)]$-coaction \eqref{lc} to
obtain the Hopf fibration itself, we turn to the remaining
symmetry determined by the coaction \eqref{rc}. 
We aim at a quantum symmetry group $\Sp_q(2)$ which extends the coaction $\eqref{rc}$ on the two copies 
$\A[\HH_q]$ and $\A[\HH_q']$ in $\A[\HH^2_q]$ while being compatible with the quaternionic structure in \eqref{eqn J}.
Motivated by a general strategy ({\em cf}. \cite{wang, sol} and in
particular \cite{lprs:ncfi} for a similar case), we first define a bialgebra $\A[\M_q(2,\HH)]$ to be the
universal bialgebra for which $\A[\HH^2_q]$ is a comodule
$*$-algebra, defined as follows.

\begin{defn} 
We say that a $*$-bialgebra
$\mathcal{B}$ is a {\em transformation bialgebra} for $\A[\HH^2_q]$
if there is a $*$-algebra map
$$
\Delta_L:\A[\HH^2_q]\to \mathcal{B}\otimes\A[\HH^2_q]
$$
which is an intertwiner for the quaternionic structure, namely
\begin{equation}\label{J-int}
(\id\otimes\mathrm{J}_q)\circ \Delta_L=\Delta_L\circ\mathrm{J}_q.
\end{equation}
We then define $\A[\M_q(2,\HH)]$ as the universal transformation
bialgebra for $\A[\HH^2_q]$, in the sense that whenever
$\mathcal{B}$ is a transformation bialgebra for $\A[\HH^2_q]$ there
is a morphism of transformation bialgebras ({\em i.e.} commuting
with the coactions) from $\A[\M_q(2,\HH)]$ onto $\mathcal{B}$.
\end{defn}

From the requirement that $\A[\HH^2_q]$ be a comodule algebra, we
derive the structure of $\A[\M_q(2,\HH)]$. In order to
satisfy the universality condition, we see  that
$\A[\M_q(2,\HH)]$ is generated as a $*$-algebra by the
entries of a $4\times 4$ matrix $A=(A_i{}^j)$. Then in terms of the
matrix \eqref{u matrix} (for the time being forgetting the sphere relations)
there is a linear map
\begin{equation}\label{unico}
\Delta_L:\A[\HH^2_q]\to\A[\M_q(2,\HH)]\otimes\A[\HH^2_q],\qquad
{\su}_{i}{}^a\mapsto (A_\mu{}^i)^*{}  \otimes {\su}_{\mu}{}^a ,
\end{equation}
for $i=1,\ldots,4$ and $a=1,2$, extended as an algebra map. This
 becomes a left $\A[\M_q(2,\HH)]$-coaction if we equip
$\A[\M_q(2,\HH)]$ with the (co-opposite) matrix coalgebra structure 
$$
\Delta^\cop(A_i{}^j)=A_\mu{}^j\otimes A_i{}^\mu, \qquad
\ep(A_i{}^j)=\delta_i^j ,
$$
on generators for $i,j=1,\dots,4$, each extended as $*$-algebra maps.

Now, imposing the condition \eqref{J-int} shows that $A$ necessarily
has the form
\begin{equation}\label{A form}
A=\begin{pmatrix}\sa_i{}^j&\sb_i{}^{j}\\\sc_{i}{}^j&\sd_{i}{}^{j}\end{pmatrix}=\begin{pmatrix}a_1&-q a_2^*&b_1&-q b_2^*\\a_2&a_1^*&b_2&b_1^*\\c_1&-q c_2^*&d_1&-q d_2^*\\c_2&c_1^*&d_2&d_1^*\end{pmatrix}
\end{equation}
so that, in some sense, $A$ may be thought of as a deformed $2\times
2$ matrix of quaternions,
$$
A=\begin{pmatrix}\bfa&\mathbf{b}\\\mathbf{c}&\mathbf{d}\end{pmatrix},
\qquad
\text{where}~~\bfa=(\sa_i{}^j)=\begin{pmatrix}a_1&-q a_2^*\\a_2&a_1^*\end{pmatrix},
$$
with similar notation for the remaining blocks
$\mathbf{b},\mathbf{c},\mathbf{d}$. The requirement that $\Delta_L$
be an algebra map allows us to deduce the algebra structure of
$\A[\M_q(2,\HH)]$.

\begin{prop}\label{pr:algrels}
In terms of the $R$-matrix \eqref{R mat}, with notation $\zeta=q^{-1/2}$, the
algebra relations in $\A[\M_q(2,\HH)]$ are given by
\begin{eqnarray*}
R_i{}^\alpha{}_k{}^\beta  \sa_\alpha{}^j\sa_\beta{}^l=
 \sa_k{}^\beta\sa_i{}^\alpha R_\alpha{}^j{}_\beta{}^l, & \qquad
R_i{}^\alpha{}_k{}^\beta\sb_\alpha{}^{j}\sb_\beta{}^{l}=
\sb_k{}^{\beta}\sb_i{}^{\alpha} R_{\alpha}{}^j{}_{\beta}{}^l, \\
\sb_i{}^1\sa_k{}^1=\zeta R_k{}^\alpha{}_i{}^\beta \sa_\alpha{}^1\sb_\beta{}^1,& \qquad \sb_i{}^1\sa_k{}^2=\zeta R_k{}^\alpha{}_i{}^\beta \sa_\alpha{}^2\sb_\beta{}^1, \\
\sa_i{}^1\sb_k{}^2=\zeta R_k{}^\alpha{}_i{}^\beta \sb_\alpha{}^2\sa_\beta{}^1,& \qquad \sa_i{}^2\sb_k{}^2=\zeta R_k{}^\alpha{}_i{}^\beta \sb_\alpha{}^2\sa_\beta{}^2,
\end{eqnarray*}
\begin{eqnarray*}
\sa_1{}^i\sc_1{}^k=\zeta \sc_1{}^\beta\sa_1{}^\alpha R_\alpha{}^i{}_\beta{}^k ,& \qquad 
\sa_1{}^i\sc_2{}^k=\zeta \sc_2{}^\beta\sa_1{}^\alpha R_\alpha{}^i{}_\beta{}^k, \\
\sc_1{}^i\sa_2{}^k=\zeta \sa_2{}^\beta\sc_1{}^\alpha R_\alpha{}^i{}_\beta{}^k ,& \qquad 
\sc_2{}^i\sa_2{}^k=\zeta \sa_2{}^\beta\sc_2{}^\alpha R_\alpha{}^i{}_\beta{}^k, &\qquad \\
\sb_1{}^i\sd_1{}^k=\zeta \sd_1{}^\beta\sb_1{}^\alpha R_\alpha{}^i{}_\beta{}^k ,& \qquad 
\sb_1{}^i\sd_2{}^k=\zeta \sd_2{}^\beta\sb_1{}^\alpha R_\alpha{}^i{}_\beta{}^k, \\
\sd_1{}^i\sb_2{}^k=\zeta \sb_2{}^\beta\sd_1{}^\alpha R_\alpha{}^i{}_\beta{}^k ,& \qquad 
\sd_2{}^i\sb_2{}^k=\zeta \sb_2{}^\beta\sd_2{}^\alpha R_\alpha{}^i{}_\beta{}^k, 
\end{eqnarray*}
\begin{eqnarray*}
R_i{}^\alpha{}_k{}^\beta  \sc_\alpha{}^j\sc_\beta{}^l=
\sc_k{}^\beta\sc_i{}^\alpha R_\alpha{}^j{}_\beta{}^l, & \qquad
R_i{}^\alpha{}_k{}^\beta\sd_\alpha{}^{j}\sd_\beta{}^{l}=
\sd_k{}^{\beta}\sd_i{}^{\alpha} R_{\alpha}{}^j{}_{\beta}{}^l, \\
\sd_i{}^1\sc_k{}^1=\zeta R_k{}^\alpha{}_i{}^\beta \sc_\alpha{}^1\sd_\beta{}^1, & \qquad \sd_i{}^1\sc_k{}^2=\zeta R_k{}^\alpha{}_i{}^\beta \sc_\alpha{}^2\sd_\beta{}^1, \\
\sc_i{}^1\sd_k{}^2=\zeta R_k{}^\alpha{}_i{}^\beta \sd_\alpha{}^2\sc_\beta{}^1, & \qquad \sc_i{}^2\sd_k{}^2=\zeta R_k{}^\alpha{}_i{}^\beta \sd_\alpha{}^2\sc_\beta{}^2,
\end{eqnarray*}
for $i,j,k,l=1,2$, together with relations, for all $i,j=1,2$, given by
\begin{align*}
\sc_1{}^j \sb_i{}^1&=\zeta^{2}\sb_i{}^1\sc_1{}^j, & \sd_1{}^1\sa_i{}^j-\sa_i{}^j\sd_1{}^1&=\zeta(q-\iq)\sb_i{}^1\sc{}_1{}^j,\\
\sc_2{}^j\sb_i{}^2&=\zeta^{-2}\sb_i{}^2\sc_2{}^j, &  \sa_i{}^j\sd_2{}^2-\sd_2{}^2\sa_i{}^j&=\zeta(q-\iq)\sc_2{}^j\sb{}_i{}^2,\\
\sd_1{}^2\sa_i{}^j&=\zeta^{2}\sa_i{}^j\sd_1{}^2, & \sc_1{}^j\sb_i{}^2-\sb_i{}^2\sc_1{}^j&=\zeta(q-\iq)\sa_i{}^j\sd_1{}^2,\\
\sd_2{}^1\sa_i{}^j&=\zeta^{-2}\sa_i{}^j\sd_2{}^1, & \sc_2{}^j\sb_i{}^1-\sb_i{}^1\sc_2{}^j&=\zeta(q-\iq)\sd_2{}^1\sa_i{}^j \ .
\end{align*} 
\end{prop}

\proof The left coaction \eqref{unico} is expressed in terms of
the $2\times 2$ blocks $\bfx$, $\bfx'$ as
$$
\sx_i{}^j\mapsto (\sa_\mu{}^i)^* \otimes\sx_\mu{}^j+ (\sb_\mu{}^i)^* \otimes\sx'_\mu{}^j,\qquad
\sx'_i{}^j\mapsto (\sc_\mu{}^i)^*\otimes\sx_{\mu}{}^j+ (\sd_\mu{}^i)^* \otimes\sx'_\mu{}^j
$$
for $i,j=1,2$. Applying this coaction to the relations \eqref{quat rels2} we find that, in order for the coaction $\Delta_L$ to be an algebra map, we must have
$$
R_k{}^\beta{}_i{}^\alpha\Delta_L(\sx_\alpha{}^j\sx_\beta{}^l)=\Delta_L(\sx_k{}^\beta\sx_i{}^\alpha)R_\alpha{}^j{}_\beta{}^l
$$
for all $i,j,k,l=1,2$. Expanding this in terms of the $2\times 2$ blocks $\bfx$, $\bfx'$ as above gives
\begin{multline*}
 R_k{}^\beta{}_i{}^\alpha\left((\sa_\nu{}^\beta\sa_\mu{}^\alpha)^*\otimes\sx_\mu{}^j\sx_\nu{}^l
+(\sb_\nu{}^\beta\sa_\mu{}^\alpha)^*\otimes\sx_\mu{}^j{\sx}'_\nu{}^l+\right.\\
 \qquad\qquad\qquad\qquad\left.+(\sa_\nu{}^\beta\sb_\mu{}^\alpha)^*\otimes\sx'_\mu{}^j{\sx}_\nu{}^l+(\sb_\nu{}^\beta\sb_\mu{}^\alpha)^*\otimes\sx'_\mu{}^j{\sx}'_\nu{}^l\right)
\\  = \left((\sa_\beta{}^i\sa_\alpha{}^k)^*\otimes\sx_\alpha{}^\nu\sx_\beta{}^\mu+(\sb_\beta{}^i\sa_\alpha{}^k)^*\otimes\sx_\alpha{}^\nu{\sx}'_\beta{}^\mu+\right.\\
 \qquad\qquad\qquad\qquad\left.+(\sa_\beta{}^i\sb_\alpha{}^k)^*\otimes\sx'_\alpha{}^\nu{\sx}_\beta{}^\mu+(\sb_\beta{}^i\sb_\alpha{}^k)^*\otimes\sx'_\alpha{}^\nu{\sx}'_\beta{}^\mu\right)R_\mu{}^j{}_\nu{}^l.
\end{multline*}
Since products of the form $\sx_i{}^j\sx_k{}^l$, $\sx_i{}^j\sx'_k{}^l$ and $\sx'_i{}^j\sx'_k{}^l$ are linearly independent from one another, the latter condition reduces to the simultaneous equations
\begin{align}
&R_k{}^\beta{}_i{}^\alpha(\sa_\nu{}^\beta\sa_\mu{}^\alpha)^*\otimes\sx_\mu{}^j\sx_\nu{}^l=(\sa_\beta{}^i\sa_\alpha{}^k)^*\otimes\sx_\alpha{}^\nu\sx_\beta{}^\mu R_\mu{}^j{}_\nu{}^l, \label{r1}
\end{align}
\begin{align}
& R_k{}^\beta{}_i{}^\alpha\left((\sb_\nu{}^\beta\sa_\mu{}^\alpha)^*\otimes\sx_\mu{}^j{\sx}'_\nu{}^l+(\sa_\nu{}^\beta\sb_\mu{}^\alpha)^*\otimes\sx'_\mu{}^j{\sx}_\nu{}^l\right)\label{r2}\\
&\qquad\qquad =\left((\sb_\beta{}^i\sa_\alpha{}^k)^*\otimes\sx_\alpha{}^\nu{\sx}'_\beta{}^\mu+(\sa_\beta{}^i\sb_\alpha{}^k)^*\otimes\sx'_\alpha{}^\nu{\sx}_\beta{}^\mu\right)R_\mu{}^j{}_\nu{}^l, \nn 
\end{align}
\begin{align}
& R_k{}^\beta{}_i{}^\alpha
(\sb_\nu{}^\beta\sb_\mu{}^\alpha)^*\otimes\sx'_\mu{}^j{\sx}'_\nu{}^l=(\sb_\beta{}^i\sb_\alpha{}^k)^*\otimes\sx'_\alpha{}^\nu{\sx}'_\beta{}^\mu R_\mu{}^j{}_\nu{}^l \label{r3} ,
\end{align}
for all $i,k,j,l=1,2$. Applying the relations \eqref{quat rels2} to eq.~\eqref{r1} and taking conjugates gives
$$
R_k{}^\beta{}_i{}^\alpha\sa_\nu{}^\beta\sa_\mu{}^\alpha\otimes\sx_\mu{}^j\sx_\nu{}^l=\sa_\beta{}^i\sa_\alpha{}^k\otimes\sx_\mu{}^j\sx_\nu{}^l R_\alpha{}^\nu{}_\beta{}^\mu.
$$
Next we use the fact that $R_k{}^\beta{}_i{}^\alpha=R_\alpha{}^i{}_\beta{}^k$ for all $i,k,\alpha,\beta=1,2$ and observe that, for $j\leq l$, the generators $\sx_i{}^j\sx_k{}^l$ can be taken
to be all linearly independent. As a consequence, for each fixed $\mu$, $\nu$, we must have
$$
R_\mu{}^\beta{}_\nu{}^\alpha\sa_\beta{}^i\sa_\alpha{}^k=\sa_\nu{}^\beta\sa_\mu{}^\alpha R_\alpha{}^i{}_\beta{}^k,
$$
which are just the $\bfa$-$\bfa$ relations as stated in the proposition. The $\bfb$-$\bfb$ relations are obtained in the same way from eq.~\eqref{r3}. Similarly, applying the relations \eqref{quat rels4} to eq.~\eqref{r2} and then using the fact that the generators $\sx_i{}^j\sx'_k{}^l$ for $j\leq l$ may be taken to be all linearly independent yields the $\bfa$-$\bfb$ relations as stated.
The remaining relations are obtained with the same strategy to
products of elements of the form $\sx_{i}{}^j\sx'_k{}^l$,
$\sx'_i{}^j\sx_{k}{}^l$ and $\sx'_{i}{}^j\sx'_{k}{}^l$.\endproof

It is not difficult to see
that $\A[\M_q(2,\HH)]$ is indeed the universal transformation
bialgebra for $\A[\HH^2_q]$, since the commutation relations in
Prop.~\ref{pr:algrels} and the $*$-structure \eqref{A form}
are derived from the minimal requirements that $\Delta_L$ be a
$*$-algebra map such that the compatibility \eqref{J-int} with the quaternionic structure holds.

\begin{rem}\textup{
We observe for later use that the two $*$-subalgebras generated by the $2\times 2$ blocks $\mathbf{a}$ and $\mathbf{d}$ are each isomorphic to a copy of the algebra $\tH=\A[\SU^\cop_q(2)]$  in \S\ref{se:su2}.
}
\end{rem}

In order to obtain the Hopf algebra $\A[\Sp_q(2)]$ of quantum symmetries, we quotient
$\A[\M_q(2,\HH)]$ by the two-sided $*$-ideal $I$ generated by
elements of the form
\begin{equation}\label{quot rels}
 (A_\mu{}^j)^* A_\mu{}^l-\delta_j{}^l, \qquad
  A_j{}^\mu (A_l{}^\mu)^*-\delta_j{}^l,
\end{equation}
where $j,l=1,\ldots,4$. We denote the resulting quotient algebra by
$\A[\Sp_q(2)]$ and define
$$
S:\A[\Sp_q(2)]\to \A[\Sp_q(2)],\qquad S(A_j{}^l):=(A_l{}^j)^*,
$$
extended as an anti-$*$-algebra map.

\begin{prop}The datum $(\A[\Sp_q(2)],\Delta^\cop,\ep,S)$ constitutes a
Hopf algebra.\end{prop}

\proof It is straightforward to check that 
$\Delta^\cop(I)\subset \A[\M_q(2,\HH)]\otimes I + I\otimes
\A[\M_q(2,\HH)]$ and $\ep(I)=0$, whence $I$ is a Hopf $*$-ideal and so the quotient
$\A[\Sp_q(2)]$ is a bialgebra. The very form of the ideal $I$ means
that the map $S$ satisfies the properties of an antipode.\endproof

\begin{rem}\textup{We stress that our quantum group $\Sp_q(2)$ is {\em not}
the `FRT' quantum group coming from the C-series of Lie groups
\cite{rtf:qlg}, which is used in particular in \cite{lpr:qsb} to
construct a quantum Hopf fibration. It {\em is} however a
deformation of the classical Lie group $\Sp(2)$, the symmetry group
of the classical Hopf fibration $S^7\to S^4$, our notation reflecting this.}\end{rem}

\subsection{Quantum homogeneous spaces}\label{hom spaces}
Having constructed the quantum group $\Sp_q(2)$ of symmetries of the space
$\HH^2_q$, we now have to check that its action descends to the
spheres $S^7_q$ and $S^4_q$ as a group of symmetries of the Hopf fibration.

\begin{lem} The coaction $\Delta_L:\A[\HH^2_q]\to\A[\Sp_q(2)]\otimes\A[\HH^2_q]$ preserves the two-sided $*$-ideal generated by
$r^2-1$.\end{lem}

\proof We observe that
$r^2=\sum\nolimits_\mu(\su_\mu{}^a)^*\su_\mu{}^a$ for both $a=1,2$. 
Using this, we compute that
\begin{align*}
\Delta_L(r^2)&=\sum\nolimits_\mu\Delta_L(({\su}_\mu{}^a)^*)\Delta_L({\su}_\mu{}^a)=
\sum\nolimits_{\mu,\alpha,\beta}S(A_\alpha{}^\mu)A_\mu{}^\beta\otimes ({\su}_\alpha{}^a)^*{\su}_\beta{}^a\\&=
\sum\nolimits_{\alpha,\beta}\delta_\alpha{}^\beta\otimes ({\su}_\alpha{}^a)^*{\su}_\beta{}^a= 1\otimes r^2,\end{align*} 
having used the defining relations \eqref{quot rels} for
$\A[\Sp_q(2)]$. The result now follows easily.\endproof

It follows that $\A[S^7_q]$ is a left comodule $*$-algebra for
the coaction of $\A[\Sp_q(2)]$.  By construction, this coaction
commutes with the right coaction of $\A[\SU_q(2)]$ in
eq.~\eqref{su coact}.
The fact that they commute is due to $\A[\Sp_q(2)]$ coacting upon the rows of the matrix \eqref{u matrix}, whereas $\A[\SU_q(2)]$ coacts upon the columns. In the classical case, this is nothing other than the statement that left and right matrix multiplication are mutually commuting operations.  The subalgebra $\A[S^4_q]$ is therefore
an $\A[\Sp_q(2)]$-comodule $*$-algebra as well. 

\begin{lem}
Let $\mathcal{I}$ be the two-sided $*$-ideal of $\A[\Sp_q(2)]$
generated by the elements
$$b_1,~b_1^*,~b_2,~b_2^*,~c_1,~c_1^*,~c_2,~c_2^*,~d_2,~d_2^*,~d_1-1,~d_1^*-1.$$
Then $\mathcal{I}$ is a Hopf $*$-ideal, namely
\begin{equation}
\label{hopf-id} \ep(\mathcal{I})=0,\qquad \Delta^\cop(\mathcal{I})\subset
\A[\Sp_q(2)]\otimes\mathcal{I}+\mathcal{I}\otimes\A[\Sp_q(2)],\qquad
S(\mathcal{I})\subset\mathcal{I}.
\end{equation}
\end{lem}

\proof The fact that $\mathcal{I}$ is indeed a two-sided $*$-ideal follows by inspection of the algebra relations in Prop.~\ref{pr:algrels}. The properties \eqref{hopf-id} are all easy to verify by direct
computations which we omit for the sake of brevity.\endproof

It follows that we can form the quotient Hopf algebra
$\A[\Sp_q(2)]/\mathcal{I}$. If we write
$\pi_\mathcal{I}$ for the canonical projection, given on generators by
$$
\pi_\mathcal{I}:\begin{pmatrix}\bfa&\bfb\\\bfc&\bfd\end{pmatrix}
\mapsto \begin{pmatrix}\bfa&0\\0&\mathbbm{1}_2\end{pmatrix},
$$
then the quotient may be identified with the subalgebra of $\A[\Sp_q(2)]$ generated by
$\pi_\mathcal{I}(A)=\textrm{diag}\,(\bfa,\mathbbm{1}_2)$, subject to
the relations $\bfa\bfa^*=\bfa^*\bfa=\mathbbm{1}_2$, so that the quotient is nothing other than a copy of the Hopf algebra $\A[\SU^\cop_q(2)]$.

There is a corresponding right coaction given by coproduct followed
by projection:
\begin{equation}\label{co-opp}
\Delta_\mathcal{I}:\A[\Sp_q(2)]\to\A[\Sp_q(2)]\otimes \A[\SU^\cop_q(2)],\qquad
\Delta_\mathcal{I}:=(\id\otimes\pi_{\mathcal{I}})\circ \Delta^\cop.
\end{equation}

\begin{prop}\label{isom1}There is a $*$-algebra isomorphism
$$
\phi_\mathcal{I}: \A[S^7_q]\to \A[\Sp_q(2)]^{\A[\SU^\cop_q(2)]}
$$
between $\A[S^7_q]$ and the algebra of coinvariants under the right
coaction $\Delta_\mathcal{I}$.\end{prop}

\proof Given a classical point of $\A[S^7_q]$, {\em i.e.} a
$*$-algebra map $\phi:\A[S^7_q]\to\C$, the stated isomorphism is
given by evaluating $\phi$ against the coaction
$\Delta_L:\A[S^7_q]\to\A[\Sp_q(2)]\otimes\A[S^7_q]$,
$$
\phi_\mathcal{I}:=(\id\otimes\phi)\circ\Delta_L.
$$
In particular, we choose the classical point described in
Remark~\ref{rem:cp}. On generators, the resulting isomorphism is
computed to be 
\begin{equation}\label{lin iso}
\begin{pmatrix}\iq z_1&z_2&\iq z_3&z_4\\-z_2^*&z_1^*&-z_4^*&z_3^*\end{pmatrix}^{\textup{tr}}\mapsto
\begin{pmatrix} -\iq\sc_1{}^2&\sc_1{}^1&-\iq\sd_1{}^2 &\sd_1{}^1\\-\sc_2{}^2&q\, \sc_2{}^1&-\sd_2{}^2&q\, \sd_2{}^1\end{pmatrix}^{\textup{tr}}
\end{equation}
which we extend as an algebra map.\endproof

\begin{rem}\textup{
Using the Hopf algebra surjection from $\A[\Sp_q(2)]$ to the copy of $\A[\SU^\cop_q(2)]$ generated by the $2\times 2$ block $\bfa$ of the matrix \eqref{A form}, Prop.~\ref{isom1} realises $S ^7_q$ as a quantum homogeneous space. 
On the other hand, by inspection of the relations in Prop.~\ref{pr:algrels} one sees that there is no such projection onto a subalgebra of $\A[\Sp_q(2)]$ generated by the $2\times 2$ blocks $\bfa$ and $\bfd$. 
This means that, in contrast to the classical case, or the $\theta$-deformed case in \cite{lprs:ncfi}, we cannot realise $S^4_q$ as a quantum homogeneous space of $\Sp_q(2)$. A similar phenomenon occurs in the case of the quantum four-sphere in \cite{lpr:qsb}.
}
\end{rem}

However, in the same way as for the sphere $S^7_q$, we deduce the
homogeneous space structure of quantum twistor space $\CP^3_q$. Let
$\mathcal{K}$ be the two-sided $*$-ideal of $\A[\Sp_q(2)]$ generated
by the elements $\sb_i{}^j,\sc_i{}^j,d_2,d_2^*$. Just as above it
follows that $\mathcal{K}$ is a Hopf $*$-ideal. We write
$\pi_\mathcal{K}:\A[\Sp_q(2)]\to\A[\Sp_q(2)]/\mathcal{K}$ for the
canonical projection to the quotient.  It is clear that
the quotient is generated as a $*$-algebra by the
entries of the matrix
$$
\pi_{\mathcal{K}}(A)=\begin{pmatrix}a_1 & -q a_2^* & 0 & 0 \\ a_2 & a_1^* & 0 & 0 \\ 0 & 0 & d_1 & 0 \\ 0 & 0 & 0 & d_1^*\end{pmatrix}
$$ 
subject to the relations
$\bfa^*\bfa=\bfa\bfa^*=\mathbbm{1}_2$, $d_1d_1^*=d_1^*d_1=1$ and it is
hence isomorphic to the Hopf algebra $\A[\SU^\cop_q(2)]\otimes\A[\U(1)]$.
Here there is a right coaction:
\begin{equation}\label{K-co}
\Delta_\mathcal{K}:\A[\Sp_q(2)]\to\A[\Sp_q(2)]\otimes \left(\A[\SU^\cop_q(2)]\otimes\A[\U(1)]\right), \quad
\Delta_\mathcal{K}:=(\id\otimes\pi_{\mathcal{K}})\circ \Delta^{\textup{cop}}.
\end{equation}
This construction results in the following.

\begin{prop}\label{twistor hom}There is a $*$-algebra isomorphism
$$
\phi_\mathcal{K}:\A[\CP^3_q]\to\A[\Sp_q(2)]^{\A[\SU^\cop_q(2)]\otimes\A[\U(1)]}
$$
between $\A[\CP^3_q]$ and the algebra of coinvariants under the right
coaction \eqref{K-co}.
\end{prop}

\proof From the above it is clear that the coinvariants under
$\Delta_\mathcal{K}$ are precisely the $\U(1)$-invariant elements in
the last two columns of the matrix $A$, which in turn may be identified with the
entries of the matrix \eqref{twistor proj}, yielding
\begin{align}\label{twistor
ident}
\qp=\begin{pmatrix}-\iq\sc_1{}^2\sc_2{}^1&-\iq\sc_1{}^2\sc_2{}^2&-\iq\sc_1{}^2\sd_2{}^1&-\iq\sc_1{}^2\sd_2{}^2\\
\sc_1{}^1\sc_2{}^1&\sc_1{}^1\sc_2{}^2&\sc_1{}^1\sd_2{}^1&\sc_1{}^1\sd_2{}^2\\
-\iq\sd_1{}^2\sc_2{}^1&\iq\sd_1{}^2\sc_2{}^2&-\iq\sd_1{}^2\sd_2{}^1&\iq\sd_1{}^2\sd_2{}^2\\
\sd_1{}^1\sc_2{}^1&\sd_1{}^1\sc_2{}^2&\sd_1{}^1\sd_2{}^1&\sd_1{}^1\sd_2{}^2\end{pmatrix}.\end{align}
This expression will be useful when we consider
differential structure on twistor space.\endproof

\subsection{A differential calculus on $\Sp_q(2)$}
In this section we give a differential calculus on the quantum
group $\A[\Sp_q(2)]$, which we shall use later on to construct differential structures on the spaces $S^7_q$, $\CP^3_q$ and $S^4_q$. 

For simplicity here we use the shorthand  $P:=\A[\Sp_q(2)]$, with $\ep_P$
denoting its counit. 
Recall from \S\ref{se:dcqg} that left covariant differential calculi on $P$ are given by right ideals 
of the augmentation ideal $P^+=\textrm{Ker}\,\ep_P$. 
Thus, to give a calculus on $P$ we simply define a right ideal $I_P$ of $P^+$ 
by equipping each of the $2\times 2$
blocks generated by $\bfa$ and $\bfd$ with a copy of the ideal
$I_{\tH}$ defined in Ex.~\ref{ex:4dopp} (and corresponding to the 4D$_+$ calculus
on $\SU^\cop_q(2)$), then straightforwardly extending this quadratically  
in the simplest way. The resulting ideal $I_P$ is then 
generated by the nine elements
\begin{align}\label{a-block}
&\sa_1{}^2(\sa_1{}^1-\sa_2{}^2),~\sa_2{}^1(\sa_1{}^1-\sa_2{}^2),~(\sa_1{}^1)^2+q^{2}(\sa_2{}^2)^2-(1+q^2)(\sa_1{}^1\sa_2{}^2+\iq
\sa_1{}^2\sa_2{}^1), \\
&(\sa_1{}^2)^2,~(\sa_2{}^1)^2,~\st_\sa\sa_1{}^2,~\st_\sa\sa_2{}^1,~\st_\sa(\sa_1{}^1-\sa_2{}^2),~\st_\sa(q^2\sa_1{}^1+\sa_2{}^2-(q^2+1)), \nonumber
\end{align}
where $\st_\sa:=q^2\sa_1{}^1+\sa_2{}^2-(q^3+q^{-1})$; the nine
elements
\begin{align}\label{d-block}
&\sd_1{}^2(\sd_1{}^1-\sd_2{}^2),~\sd_2{}^1(\sd_1{}^1-\sd_2{}^2),~(\sd_1{}^1)^2+q^{2}(\sd_2{}^2)^2-(1+q^2)(\sd_1{}^1\sd_2{}^2+\iq
\sd_1{}^2\sd_2{}^1), \\
&(\sd_1{}^2)^2,~(\sd_2{}^1)^2,~\st_\sd\sd_1{}^2,~\st_\sd\sd_2{}^1,~\st_\sd(\sd_1{}^1-\sd_2{}^2),~\st_\sd(q^2\sd_1{}^1+\sd_2{}^2-(q^2+1)),\nonumber
\end{align}
where $\st_\sd:=q^2\sd_1{}^1+\sd_2{}^2-(q^3+q^{-1})$; and the
elements
\begin{align}\label{extras}
&\sb_{i}{}^{j}\sb_{k}{}^{l},\quad\sc_{i}{}^j\sc_{k}{}^l,\quad\sb_{i}{}^{j}\sc_{k}{}^l,\quad \sb_{k}{}^{l}(\sa_i{}^j-\delta_i^j),\\
&\sc_{k}{}^{l}(\sa_i{}^j-\delta_i^j),\quad\sb_{k}{}^{l}(\sd_{i}{}^{j}-\delta_{i}^{j}),\quad\sc_{k}{}^{l}(\sd_{i}{}^{j}-\delta_{i}^{j}),\quad
(\sa_i{}^j-\delta_i^j)(\sd_k{}^l-\delta_k^l), \nonumber
\end{align}
for $i,j,k,l=1,2$. 

\begin{defn}\label{def:qgcalc}
We write $\Omega^1\Sp_q(2)$ for the left-covariant differential
calculus on the Hopf algebra $P=\A[\Sp_q(2)]$ corresponding to the
right ideal $I_P$.
\end{defn}

Let us examine the calculus $\Omega^1\Sp_q(2)$ more closely by looking at the space of corresponding 
left invariant one-forms given in the usual way \cite{wor}. 
In parallel with eqs.~\eqref{4dspan} for the calculus $\Omega^1\SU^\cop_q(2)$, let $\sa_0$, $\sa_z$, $\sd_0$, $\sd_z$ be the elements defined by
the equations
\begin{align}\label{az1}
\sa_1{}^1-1&=(\iq-1)\sa_0+(q-1)\sa_z, &
\sa_2{}^2-1&=(q-1)\sa_0+(\iq-1+\sigma^2\iq)\sa_z, \nonumber \\
\sd_1{}^1-1&=(\iq-1)\sd_0+(q-1)\sd_z, &
\sd_2{}^2-1&=(q-1)\sd_0+(\iq-1+\sigma^2\iq)\sd_z.
\end{align}

\begin{lem}
The vector space $P^+/I_P$ has a basis given by the
elements $\sa_1{}^2$, $\sa_0$, $\sa_2{}^2$, $\sa_z$, $\sd_1{}^2$,
$\sd_0$, $\sd_2{}^2$, $\sd_z$ and the elements $\sc_i{}^j$ for $i,j=1,2$.
\end{lem}
\proof We define an equivalence relation $\sim$ on $P^+$
by setting $x\sim y$ if and only if $x-y\in I_P$. Then from the
form of the generators \eqref{extras} and the
eqs.~\eqref{az1}, we immediately see that
\begin{align*}
\sa_1{}^2\sb_k{}^l\sim \sa_0\sb_k{}^l\sim \sa_2{}^1\sb_k{}^l\sim
\sa_z\sb_k{}^l&\sim 0, & \sd_1{}^2\sb_k{}^l\sim \sd_0\sb_k{}^l\sim
\sd_2{}^1\sb_k{}^l\sim \sd_z\sb_k{}^l&\sim 0,
\\
\sa_1{}^2\sc_k{}^l\sim \sa_0\sc_k{}^l\sim \sa_2{}^1\sc_k{}^l\sim
\sa_z\sc_k{}^l&\sim 0, & \sd_1{}^2\sc_k{}^l\sim \sd_0\sc_k{}^l\sim
\sd_2{}^1\sc_k{}^l\sim \sd_z\sc_k{}^l&\sim 0,
\end{align*}
together with $\sb_i{}^j\sb_k{}^l\sim \sb_i{}^j\sc_k{}^l\sim
\sc_i{}^j\sc_k{}^l\sim 0$ for all $i,j,k,l=1,2$. 

It follows that any
expression which is quadratic in elements from the off-diagonal
blocks $\bfb$, $\bfc$ can be rewritten as expressions which are
linear in the generators in $\bfb$ and $\bfc$. Moreover, we know
from \S\ref{se:su2} 
that all expressions which are quadratic in the
elements $\sa_0$, $\sa_z$, $\sa_1{}^2$, $\sa_2{}^1$ can be made
linear; similarly for expressions which are quadratic in the
elements $\sd_0$, $\sd_z$, $\sd_1{}^2$, $\sd_2{}^1$. The fact that
all products of the form
$(\sa_i{}^j-\delta_i^j)(\sd_k{}^l-\delta^k_l)$ are equivalent to
zero means that every quadratic polynomial can be rewritten as a
linear combination of the generators stated in the lemma. It follows
by induction that polynomials of arbitrary degree can be made
linear. 
Finally, with $P^+$ being stable under the antipode, from the relations $S(\sb_i{}^j)=(\sc_j{}^i)^*$ we deduce that the generators $\sb_i{}^j$ are not independent from the generators $\sb_i{}^j$ and so we do not need them.
\endproof

We deduce that the calculus $\Omega^1\Sp_q(2)$ is spanned as a left module by a twelve-dimensional space of left-invariant one-forms, in contrast with the classical dimension of ten. 

Not needing the commutation relations and the bimodule structure of the calculus in what follows we shall not give them here.

\section{Framed Quantum Manifolds}\label{se:non-univ}
According to the frame bundle theory of \cite{Ma1} which we sketched in \S\ref{se:frame}, a choice of differential structure on the Hopf algebra $P$ which satisfies the compatibility conditions \eqref{hom bundle} is sufficient to guarantee that the associated homogeneous space $M=P^H$ is a framed quantum manifold. In this section we shall see how this is indeed the case for the spaces $S^7_q$ and $\CP^3_q$, both of them being quantum homogeneous spaces of $\Sp_q(2)$. We then use this to equip the quantum four-sphere $S^4_q$ --- which, in contrast, is not a quantum homogeneous space --- with a finite-dimensional differential calculus.

\subsection{Framed manifold structure of $S^7_q$}
First of all we compute the framing of the quantum seven-sphere
$S^7_q$. We use the structure of $S^7_q$ described
in Prop.~\ref{isom1}, writing 
$$P=\A[\Sp_q(2)],\qquad \tH=\A[\SU^\cop_q(2)], \qquad \tM=\A[S^7_q],
$$
so that we have $\tM=P^{\tH}$ as a
quantum homogeneous space.

The Hopf algebra $P=\A[\Sp_q(2)]$ is equipped with the
left-covariant differential calculus $\Omega^1P$ determined by the
right ideal $I_P$ according to Definition~\ref{def:qgcalc}. The
fibre quantum group $\tH=\A[\SU^\cop_q(2)]$ is equipped with the
bicovariant differential calculus $\Omega^1\tH$ determined by the
ideal $I_{\tH}:=\pi_\mathcal{I}(I_P)$. By inspection, we see that the calculus $\Omega^1\tH$ is simply a copy
of the 4D calculus on the algebra $\A[\SU^\cop_q(2)]$ defined in Example~\ref{ex:4dopp}.

\begin{lem}\label{cond1}With $P$, $\tH$, $\tM$ as above, the differential calculi
$\Omega^1P$, $\Omega^1\tH$, $\Omega^1\tM$ satisfy the (non-universal calculi) compatibility
conditions of eq.~\eqref{hom bundle}.
\end{lem}

\proof The condition $\pi_\mathcal{I}(I_P)=I_{\tH}$ holds by
definition of the differential calculus on $\tH$. Since the 4D
calculus on $\SU^\cop_q(2)$ is bicovariant, its defining ideal is
$\textup{Ad}_R$-stable, thus the
generators of $I_P$ in eqs.~\eqref{a-block} and \eqref{d-block} are
$\textrm{Ad}_R$-stable. Moreover, one finds that
\begin{align*}
(\id\otimes\pi_{\mathcal{I}})\textup{Ad}_R(\sa_i{}^{j})&=\sa_\mu{}^\nu\otimes S(\sa_\nu{}^j)\sa_i{}^\mu , &
(\id\otimes\pi_{\mathcal{I}})\textup{Ad}_R(\sd_i{}^{j})&=\sd_i{}^j\otimes 1, \\
(\id\otimes\pi_{\mathcal{I}})\textup{Ad}_R(\sb_i{}^{j})&=\sb_\mu{}^j\otimes
\sa_i{}^\mu, &
(\id\otimes\pi_{\mathcal{I}})\textup{Ad}_R(\sc_{i}{}^j)&=\sc_i{}^{\nu}\otimes
S(\sa_\nu{}^j).
\end{align*}
Using these, the $\textrm{Ad}_R$-condition \eqref{hom bundle}
for the generators in eq.~\eqref{extras} is easily
verified.\endproof

This means that we may apply the framing theory of \S\ref{se:frame}
to express $S^7_q$ as a framed quantum manifold, immediately yielding the following theorem.

\begin{thm} \label{th:s7}
In terms of the expressions for the framing as in eq.~\eqref{frame}, the quantum homogeneous 
space $\A[S^7_q]=P^{\tH}$ is a framed quantum
manifold with cotangent bundle 
\begin{equation}\label{s7-iso}
\Omega^1S^7_q \simeq (P\otimes\tV)^{\tH},
\end{equation}
where the right ${\tH}$-comodule $\tV$ is determined to be
\begin{equation}\label{V7s}
\tV=\textup{Span}_{\C}\,\left\{ \sc_1{}^1, \sc_1{}^2,
\sc_2{}^1, \sc_2{}^2,\sd_1{}^1-1,\sd_1{}^2, \sd_2{}^1, \sd_2{}^2-1\right\}
\end{equation}
in terms of the generators in the $4\times 4$ matrix \eqref{A form}.
\end{thm}

\proof The framing comodule $\tV$ is computed using eq.~\eqref{frame} as follows. From Prop.~\ref{isom1},
it is clear that $P^+\cap \tM$ is just $\tM^+:=\textup{Ker}\,\ep_{\tM}$, the
restriction of the counit $\ep_P$ to $\tM$. In our case, with $\tM$
being the algebra generated by the last two rows of 
\eqref{A form}, we have
$$
\tM^+= \la\sc_1{}^1, \sc_1{}^2,
\sc_2{}^1, \sc_2{}^2,\sd_1{}^2, \sd_0, \sd_2{}^1, \sd_z\ra
$$
as a right ideal. The ideal $I_P\cap \tM$ is generated by the
elements in \eqref{d-block} together with the elements in
\eqref{extras} of the form $\sc_i{}^j\sc_k{}^l$,
$\sc_k{}^l(\sd_i{}^j-\delta_i{}^j)$. By inspection we deduce that
$I_P\cap \tM$ contains all elements which are quadratic in the
generators of $\tM^+$. It follows that the comodule $\tV$ is the
eight-dimensional vector space given above.

The right $\A[\SU^\cop_q(2)]$-coaction on $\tV$ is computed from the formula
in \eqref{frame} to be
\begin{equation*}
\Delta_R(\sd_i{}^j-\delta_i{}^j)=(\sd_i{}^j-\delta_i{}^j)\otimes 1,\qquad
\Delta_R(\sc_{i}{}^j)=\sc_i{}^{\mu}\otimes S(\sa_{\mu}{}^j), \qquad i,j=1,2.
\end{equation*}
The corresponding soldering form is computed using the formula in eq.~\eqref{frame} to be
\begin{align*} 
\theta(\sc_i{}^j)&=S(\sa_\mu{}^j)\D\sc_i{}^\mu+S(\sc_\mu{}^j)\D\sd_i{}^\mu, \\
\theta(\sd_i{}^j-\delta_i^j)&=
S(\sb_\mu{}^j)\D\sc_i{}^\mu+S(\sd_\mu{}^j)\D\sd_i{}^\mu,
\end{align*}
for each $i,j=1,2$.\endproof

In this way, we equip the seven-sphere $S^7_q$ with
an eight-dimensional differential calculus. The increase in
dimension from the classical value is the price we have to pay for
having a bicovariant differential calculus on the structure
group $\A[\SU^\cop_q(2)]$ of the fibration.

It is useful to see how the geometry of the calculus $\Omega^1S^7_q$ is reflected in the frame bundle construction of 
Theorem~\ref{th:s7}. Due to its proof, the calculus is spanned, as a left $\A[S^7_q]$-module, by the following $\tH$-invariant elements of the tensor product comodule $P\otimes\tV$:
\begin{align}\label{invs}
\sa_\mu{}^j\otimes\sc_i{}^\mu,\qquad \sb_\mu{}^j\otimes\sc_i{}^\mu, \qquad 1\otimes(\sd_i{}^j-\delta_i{}^j),\qquad i,j=1,2.
\end{align}
Under the isomorphism \eqref{s7-iso} induced by the map $\varpi$ of eq.~\eqref{iso-cl}, these are respectively carried onto the following elements of $\Omega^1S^7_q$:
\begin{align}\label{comps}
\sa_\mu{}^j\otimes\sc_i{}^\mu&~\mapsto~ \sa_\mu{}^j\left(S(\sa_\nu{}^\mu)\otimes \sc_i{}^\nu+S(\sc_\nu{}^\mu)\otimes \sd_i{}^\nu\right),\nonumber\\
\sb_\mu{}^j\otimes\sc_i{}^\mu&~\mapsto~\sc_\mu{}^j\left(S(\sa_\nu{}^\mu)\otimes \sc_i{}^\nu+S(\sc_\nu{}^\mu)\otimes \sd_i{}^\nu\right),\\
1\otimes(\sd_i{}^j-\delta_i{}^j)&~\mapsto~ S(\sb_\mu{}^j)\otimes \sc_i{}^\mu+S(\sd_\mu{}^j)\otimes\sd_i{}^\mu -1\otimes\delta_i{}^j, \qquad \textup{for} \quad i,j=1,2.
\nonumber
\end{align}
One immediately sees that, {\em a priori}, these are not elements of the kernel of the multiplication on $\A[S^7_q]$, as one would expect of the differential one-forms $\Omega^1S^7_q$. This apparent problem is resolved by taking linear combinations of these elements and using the defining properties of the Hopf algebra $\A[\Sp_q(2)]$. For example, under the mapping $\varpi$ we have
\begin{align*}
& \varpi\left(\sa_\mu{}^j\otimes\sc_i{}^\mu +\sc_\mu{}^j\otimes(\sd_i{}^\mu-\delta_\mu{}^j)\right)  \\ 
& \qquad = \sa_\mu{}^j\left(S(\sa_\nu{}^\mu)\otimes \sc_i{}^\nu+S(\sc_\nu{}^\mu)\otimes \sd_i{}^\nu\right)
 + \sc_\mu{}^j\left(S(\sb_\nu{}^\mu)\otimes \sc_i{}^\nu+S(\sd_\nu{}^\mu)\otimes\sd_i{}^\nu -1\otimes\delta_i{}^\mu\right) \\
& \qquad = \left(\sa_\mu{}^j S(\sa_\nu{}^\mu)+\sc_\mu{}^j S(\sb_\nu{}^\mu)\right)\otimes\sc_i{}^\nu+ \left(\sa_\mu{}^j S(\sc_\nu{}^\mu)+ \sc_\mu{}^j S(\sd_\nu{}^\mu)\right)\otimes \sd_i{}^\nu -\sc_\mu{}^j\otimes\delta_i{}^\mu \\
& \qquad =1\otimes\sc_i{}^j-\sc_i{}^j\otimes 1,
\end{align*}
with a similar computation yielding that 
$$
\varpi\left(\sb_\mu{}^j\otimes\sc_i{}^\mu +\sd_\mu{}^j\otimes(\sd_i{}^\mu-\delta_\mu{}^j)\right)= 1\otimes\sd_i{}^j-\sd_i{}^j\otimes 1
$$ 
for each $i,j=1,2$. These and other similar expressions will prove useful when we come to equip the four-sphere $S^4_q$ with a first order differential calculus.

\subsection{Framing of noncommutative twistor space} It is now
relatively straightforward to convert the differential calculus on
$S^7_q$ to one on $\CP^3_q$, again using the quantum framing theory.
To this end, we use the homogeneous space description of twistor space
constructed in Prop.~\ref{twistor hom}, writing
$$
P=\A[\Sp_q(2)],\qquad \hH=\A[\SU^\cop_q(2)]\otimes\A[\U(1)],\qquad \hM=\A[\CP_q^3]
$$  
for total space, structure quantum group and base space of the principal bundle, the latter being the
subalgebra $\hM=P^{\hH}$ of invariants for the coaction
$\Delta_\mathcal{K}:P\to P\otimes \hH$ in \eqref{K-co}.

Once again, the Hopf algebra $P=\A[\Sp_q(2)]$ is equipped with the
left-covariant differential calculus $\Omega^1P$ determined by the
right ideal $I_P$ according to Definition~\ref{def:qgcalc}. This
time, the structure quantum group $\hH$ is
equipped with the bicovariant differential calculus $\Omega^1\hH$
determined by the ideal $I_{\hH}:=\pi_\mathcal{K}(I_P)$. By inspection again, we see that the
five-dimensional calculus $\Omega^1\hH$ is nothing other than the tensor product
bimodule $\Omega^1\SU^\cop_q(2)\otimes\Omega^1\U(1)$ of the differential calculi defined in \S\ref{se:su2}. 

\begin{lem}\label{cond2}With $P$, $\hH$, $\hM$ as above, the differential calculi
$\Omega^1P$, $\Omega^1\hH$, $\Omega^1\hM$ satisfy the (non-universal calculi) compatibility
conditions of eq.~\eqref{hom bundle}.\end{lem}

\proof The condition $\pi_\mathcal{K}(I_P)=I_{\hH}$ holds by
definition of the differential calculus on $\hH$. We find that, under
the map $(\id\otimes\pi_{\mathcal{K}})\textup{Ad}_R$ the
generators of $I_P$ transform according to
\begin{align*}
\sd_1{}^1&\mapsto \sd_1{}^1\otimes 1, & \sd_1{}^2&\mapsto
\sd_1{}^2\otimes t^*{}^2, & \sb_i{}^{1}&\mapsto \sb_\mu{}^1\otimes \sa_i{}^\mu\, t, & \sb_i{}^{2}&\mapsto \sb_\mu{}^2\otimes \sa_i{}^\mu\, t^*, \\
\sd_2{}^1&\mapsto \sd_2{}^1\otimes t{}^2, & \sd_2{}^2&\mapsto
\sd_2{}^2\otimes 1, & \sc_{1}{}^j&\mapsto \sc_1{}^\nu\otimes 
S(\sa_\nu{}^j)t^*, & \sc_{2}{}^j&\mapsto \sc_2{}^\nu \otimes S(\sa_\nu{}^j) t,
\end{align*}
together with $\sa_i{}^j\mapsto\sa_\mu{}^\nu\otimes  S(\sa_\nu{}^j)
\sa_i{}^\mu$. In particular, it follows that
$\sd_0\mapsto \sd_0\otimes 1$ and $\sd_z\mapsto \sd_z\otimes 1$.
Using these transformation rules, the condition \eqref{hom bundle}
is easy to verify.\endproof

As before, this means that we can express $\CP^3_q$ as a framed quantum manifold.

\begin{thm}\label{th:tw-frm}
The quantum homogeneous space $\A[\CP^3_q]=P^{\hH}$ is a framed
quantum manifold with cotangent bundle 
$$
\Omega^1\CP^3_q\simeq (P\otimes\widehat V)^{\hH},
$$
where the right $\hH$-comodule $\widehat V$
is determined to be
\begin{align*}
\widehat V&=\textup{Span}_{\C}\,\left\{
\qp_{33},\qp_{i4},\qp_{4j}~|~i,j=1,2,3
\right\}\\&=\textup{Span}_{\C}\,\left\{ \sd_1{}^2\sd_2{}^1,\,\sd_1{}^1\sd_2{}^1,\,\sd_1{}^1\sc_2{}^1,\,\sd_1{}^1\sc_2{}^2,\,
\sd_1{}^2\sd_2{}^2,\,\sc_1{}^2\sd_2{}^2,\,\sc_1{}^1\sd_2{}^2,\right\},
\end{align*}
in terms of the entries of the projection $\qp$ of eq.~\eqref{twistor ident}.
\end{thm}
\proof 
To compute the comodule $\widehat V$, we observe once again that $P^+\cap \hM$ is just the
restriction $\hM^+:=\textup{Ker}\,\ep_{\hM}$ of the counit $\ep_P$ to the
subalgebra $\hM$. In this case, $\hM=\A[\CP^3_q]$ is generated by the
entries of the projection $\qp$ of eq.~\eqref{twistor ident}, whence
$\hM^+$ is generated as a left ideal by the entries of $\qp$, but
the entry $\qp_{44}=\sd_1{}^1\sd_2{}^2$, which is not in the kernel
of $\ep_{\hM}$.

Next we need the ideal $I_P\cap \hM$. Recall that $\hM=\A[\CP^3]$ is
the degree zero subalgebra of $\A[S^7_q]$ with respect to the
$\ZZ$-grading defined in eq.~\eqref{grad}. Hence, to compute
$I_P\cap \hM$, we need to find the ideal of degree zero elements in
$I_P$. By inspecting the generators of $I_P$ in
eqs.~\eqref{a-block} and \eqref{d-block} we see that the only ones
which are of homogeneous degree with respect to the $\ZZ$-grading
are the elements $(\sd_1{}^2)^2$ and $(\sd_2{}^1)^2$; the remaining
generators are not of homogeneous degree and hence the ideal that
each of them generates has no intersection with $\hM$. Similarly, we
look at the generators in \eqref{extras}: those of the form
$\sc_k{}^l(\sd_j{}^j-1)$ are not of homogeneous degree and hence
have no intersection with $\hM$. So we are left with the generators of
the form $\sc_i{}^j\sc_{k}{}^l$, $\sd_1{}^2\sc_{k}{}^l$ and
$\sd_2{}^1\sc_{k}{}^l$.

We start with the right ideal $\la (\sd_2{}^1)^2\ra$. The
elements of degree zero here include
$$
(\sd_2{}^1)^2\left\{(\sd_1{}^2)^2,(\sd_1{}^2\sd_1{}^1),(\sd_1{}^1)^2,(\sc_1{}^1)^2,(\sc_1{}^2\sc_1{}^1),(\sc_1{}^2)^2\right\},
$$
so we see that $(\qp_{33})^2$, $\qp_{33}\qp_{43}$, $(\qp_{43})^2$, $(\qp_{23})^2$, $\qp_{13}\qp_{23}$ and $(\qp_{13})^2$ are in $I_P\cap \hM$. In fact, similar
considerations for $\la (\sd_1{}^2)^2\ra$ {\em et cetera}, show that all
quadratic combinations of entries of $\qp$ (as already said, not
including $\qp_{44}$) are in $I_P\cap \hM$. Last of all, we see that
the generators $\qp_{11}$, $\qp_{12}$, $\qp_{22}$, $\qp_{13}$ and
$\qp_{23}$ are already in the ideal $I_P\cap \hM$, as are their
conjugates under the $*$-operation. As a consequence  
$\widehat V=\textup{Span}_{\C}\,\left\{
\qp_{33},\qp_{i4},\qp_{4j}~|~i,j=1,2,3 \right\}$, as stated above.

The right $\widehat H$-comodule structure on $\widehat V$ is evaluated using the
formula \eqref{frame}, yielding
\begin{align}\label{tw-fr}
\Delta_R(\qp_{33})&=\qp_{33}\otimes 1,\qquad\Delta_R(\qp_{34})=\qp_{34}\otimes t{}^2,\qquad\Delta_R(\qp_{43})=\qp_{43}\otimes t^*{}^2,\\
\Delta_R(\qp_{14})&=\iq\sc_1{}^\mu\sd_2{}^2\otimes S(\sa_\mu{}^2)t=
\qp_{14}\otimes S(\sa_2{}^2)t -\iq\qp_{24}\otimes S(\sa_1{}^2)t, \nn\\
\Delta_R(\qp_{24})&=\sc_1{}^\mu\sd_2{}^2\otimes S(\sa_\mu{}^1) t=
-q\qp_{14}\otimes S(\sa_2{}^1)t+\qp_{24}\otimes S(\sa_1{}^1)t, \nn \\
\Delta_R(\qp_{41})&=\sd_1{}^1\sc_2{}^\mu\otimes S(\sa_\mu{}^1)t^*=
\qp_{42}\otimes S(\sa_2{}^1)t^*+\qp_{41}\otimes S(\sa_1{}^1)t^*,\nn \\
\Delta_R(\qp_{42})&=\sd_1{}^1\sc_2{}^\mu\otimes S(\sa_\mu{}^2)t^*=
\qp_{42}\otimes S(\sa_2{}^2)t^*+\qp_{41}\otimes S(\sa_1{}^2)t^*. \nn
\end{align}
The soldering form is not needed in what follows and we shall not compute it explicitly.\endproof

The notable phenomenon we find in the comodule $\widehat V$ is the presence of
$\qp_{33}$ as a non-zero representative; in the classical
case it would be zero, but here it makes $\widehat V$ into a
seven-dimensional vector space, with the extra dimension being
inherited from the extra `direction' in the calculus
$\Omega^1S^7_q$. As a consequence we see that the calculus $\Omega^1\CP^3_q$ has a direct sum structure
corresponding to the decomposition
\begin{equation}\label{tw-1}
\widehat V=\C^3\oplus \C\oplus \overline\C^3
\end{equation}
as $\A[\U(1)]$-comodules: the first summand here transforms under
$t$, the third summand transforms under $t^*$, whereas the second
summand is coinvariant. This gives a decomposition of the calculus
into irreducible components,
\begin{equation}\label{tw-2}
\Omega^1\CP^3_q=\Omega^1_+ \CP^3_q \oplus \Omega^1_0 \CP^3_q \oplus\Omega^1_- \CP^3_q.
\end{equation}
In the classical limit $q\to 1$, the components $\Omega^1_\pm$
become the holomorphic and anti-holomorphic one-forms on $\CP^3$.
The extra dimension in the calculus spanned by the one-dimensional
component $\Omega^1_0$ is a purely quantum feature which is not present
in the classical limit ({\em cf}. \cite{bl:sph} for a similar phenomenon on the quantum projective line $\CP^1_q$). 

\subsection{A differential calculus on $S^4_q$}
Finally we come to describe a first order differential calculus on the sphere
$S^4_q$. In contrast to our construction of the calculi on $S^7_q$ and $\CP^3_q$, we do not have a homogeneous space structure for $S^4_q$ at our disposal and so we cannot use the usual frame bundle construction, whence we resort to more direct methods. 

Recall the right coaction of $H=\A[\SU_q(2)]$ on $\A[S^7_q]$ defined in eq.~\eqref{su coact}, namely
\begin{equation}\label{su0}
\delta_R:\A[S^7_q]\to\A[S^7_q]\otimes H,\qquad \sx_i{}^j \mapsto\sx_i{}^\mu\otimes\sa_\mu{}^j,\qquad \sx'_k{}^l\mapsto \sx'_k{}^\mu\otimes\sa_\mu{}^l,
\end{equation}
for each $i,j,k,l=1,2$. The subalgebra of invariant elements under this coaction is the algebra $\A[S^4_q]$ of coordinate functions on the quantum four-sphere $S^4_q$. We will combine the isomorphism \eqref{s7-iso}, which realises the differential calculus $\Omega^1S^7_q$, with the coaction \eqref{su0} to obtain a differential calculus on $S^4_q$.

First of all, we need to transport the coaction \eqref{su0} along the isomorphism \eqref{lin iso}. 
A direct computation shows that for the generators $(\sc_i{}^j,\sd_k{}^l)$, the coaction \eqref{su0} reads 
\begin{equation}\label{su1}
\delta_R:\A[S^7_q]\to\A[S^7_q]\otimes H,\qquad\sc_i{}^j\mapsto \sc_\mu{}^j\otimes S(u_i{}^\mu), \qquad 
\sd_i{}^j\mapsto \sd_\mu{}^j\otimes S(u_i{}^\mu),
\end{equation}
for each $i,j=1,2$, the elements $(u_i{}^j)$ being a simple relabeling of the generators of the Hopf algebra $H=\A[\SU_q(2)]$ given in 
the defining corepresentation \eqref{m-gens}. They indeed make a unitarily equivalent corepresentation given by  
\begin{equation}\label{lin iso2}
(u_i{}^j)=\begin{pmatrix}u_1{}^1 & u_1{}^2 \\ u_2{}^1 & u_2{}^2 \end{pmatrix}:=
\begin{pmatrix}0 & 1 \\ -1 & 0 \end{pmatrix}
\begin{pmatrix} \sa_1{}^1 & \sa_1{}^2 \\ \sa_2{}^1 & \sa_2{}^2\end{pmatrix} 
\begin{pmatrix}0 & -1 \\ 1 & 0 \end{pmatrix} =
  \begin{pmatrix} \sa_2{}^2 & -\sa_2{}^1 \\ -\sa_1{}^2 & \sa_1{}^1\end{pmatrix} .
\end{equation}

The coaction \eqref{su1} extends {\em via} the tensor product coaction to the universal differential calculus $\widetilde\Omega^1S^7_q$ with the sub-bimodule of coinvariant elements being precisely the universal differential calculus $\widetilde\Omega^1S^4_q$. We shall use this fact momentarily.

The next step is to extend the right $H$-coaction \eqref{su1} from the algebra $\A[S^7_q]$ to the differential calculus $\Omega^1S^7_q$. 
However, as we saw in the discussion following Theorem~\ref{th:s7}, under the isomorphism \eqref{s7-iso}
the differential calculus $\Omega^1S^7_q$ depends also on generators of $P=\A[\Sp_q(2)]$ which do not belong to the subalgebra 
$\A[S^7_q]$. 
In order to extend the right $H$-coaction from $\A[S^7_q]$ to $(P\otimes\widetilde{V})^{\widetilde{H}}$ along the isomorphism \eqref{s7-iso}, we simply
extend the coaction \eqref{su1} from $\A[S^7_q]$ to the full algebra $\A[\Sp_q(2)]$ and work at that level. 
The natural way to do so is to define
\begin{equation}\label{sp}
\delta_R:\A[\Sp_q(2)]\to\A[\Sp_q(2)]\otimes H,\qquad \sa_i{}^j\mapsto \sa_i{}^j\otimes 1, \qquad \sb_i{}^j\mapsto \sb_i{}^j\otimes 1, 
\end{equation}
for each $i,j=1,2$, extended, together with eq.~\eqref{su1}, to products of generators in the obvious way. Moreover, since the coaction \eqref{su1} and \eqref{sp} evidently commutes with the right $\widetilde{H}$-coaction 
$\Delta_\mathcal{I}:P\to P\otimes \widetilde{H}$ in \eqref{co-opp}, it is natural to think of $P$ as a right comodule
$$
\widetilde{\delta}_R:P\to P\otimes (\widetilde{H}\otimes H).
$$
for the tensor product Hopf algebra $\widetilde{H}\otimes H$. 

We are ready to describe a first order differential calculus on the four-sphere $S^4_q$. Recall that, following the general theory of \S\ref{se:frame}, the framing comodule $\widetilde{V}$ for $\Omega^1S^7_q$ was obtained in Theorem~\ref{th:s7} as the quotient of $P^+\cap \widetilde{M}$ by the ideal $I_P\cap \widetilde{M}$. The differential calculus itself was thus obtained as the quotient of the universal differential calculus $\widetilde\Omega^1S^7_q\simeq(P\otimes P^+)^{\tH}$ by the $\tM$-sub-bimodule $(P\otimes (I_P\cap\tM))^{\tH}$. Finally, using the soldering form $\theta:V\to P\Omega^1\tM$ and the induced isomorphism 
\eqref{indiso}, the differential calculus  $\Omega^1\tM$ on $\tM=\A[S^7_q]$ was obtained as the quotient of the universal calculus $\widetilde\Omega^1S^7_q$ by the $\tM$-sub-bimodule $N_{\widetilde M}$ defined by
\begin{equation}\label{bmod}
N_{\widetilde M}=\varpi\left( (P\otimes(I_P\cap\tM))^{\widetilde H}\right),
\end{equation}
where $I_P$ is the right ideal defining the differential calculus on $P=\A[\Sp_q(2)]$ in Def.~\ref{def:qgcalc} and $\varpi:P\otimes P^+\to \widetilde\Omega^1P$ is the linear isomorphism \eqref{iso-cl}. Writing $M=\A[S^4_q]$, we define an $M$-sub-bimodule $N_M$ of the universal differential calculus $\widetilde \Omega^1 M$ by
\begin{equation}\label{bmod1}
N_M:= N_{\widetilde M}\cap\widetilde\Omega^1M
\end{equation}
and then define $\Omega^1S^4_q$ to be the quotient of the universal calculus $\widetilde \Omega^1 M$ by this sub-bimodule $N_M$. Let us check that the calculi $\Omega^1S^4_q$ and $\Omega^1S^7_q$ are compatible.

\begin{prop}
The $\tM$-bimodule $N_{\tM}$ satisfies the (quantum principal bundle) compatibility conditions \eqref{submodules}, that is to say
\begin{equation}\label{r-subs}
N_M=N_{\widetilde M}\cap\widetilde\Omega^1M,\qquad \delta_R(N_{\widetilde M})\subseteq N_{\widetilde M}\otimes H,\qquad \textup{ver}(N_{\widetilde M})=\tM\otimes I_H,
\end{equation}
where $I_H$ is the right ideal of $H^+$ defining the \textup{4D}${}_+$ differential calculus on $H=\A[\SU_q(2)]$.
\end{prop}

\proof First of all we use eq.~\eqref{bmod} to describe the bimodule $N_{\tM}$ more precisely. As already observed, the ideal $I_P\cap\tM$ is generated by the elements $\sc_i{}^j\sc_k{}^l$ and $\sc_i{}^j(\sd_k{}^l-\delta_k{}^l)$ for $i,j,k,l=1,2$, together with the elements in eq.~\eqref{d-block}. We have in particular
\begin{align}\label{con-2}
\varpi(1\otimes\sc_i{}^j\sc_k{}^l)&=S(\sa_\mu{}^l)S(\sa_\nu{}^j)\otimes \sc_i{}^\nu\sc_k{}^\mu+S(\sc_\mu{}^l)S(\sa_\nu{}^j)\otimes \sc_i{}^\nu\sd_k{}^\mu  \\
&\qquad\qquad +S(\sa_\mu{}^l)S(\sc_\nu{}^j)\otimes \sd_i{}^\nu\sc_k{}^\mu + S(\sc_\mu{}^l)S(\sc_\nu{}^j)\otimes  \sd_i{}^\nu\sd_k{}^\mu, \nonumber \\
\varpi(1\otimes\sc_i{}^j(\sd_k{}^l-\delta_k{}^l))&=S(\sb_\mu{}^l)S(\sa_\nu{}^j)\otimes \sc_i{}^\nu\sc_k{}^\mu+S(\sd_\mu{}^l)S(\sa_\nu{}^j)\otimes \sc_i{}^\nu\sd_k{}^\mu  \\
&\qquad\qquad +S(\sb_\mu{}^l)S(\sc_\nu{}^j)\otimes \sd_i{}^\nu\sc_k{}^\mu + S(\sd_\mu{}^l)S(\sc_\nu{}^j)\otimes  \sd_i{}^\nu\sd_k{}^\mu\nonumber \\
&\qquad\qquad - S(\sb_\mu{}^j)\otimes \sc_i{}^\mu -S(\sd_\mu{}^j)\otimes \sd_i{}^\mu, \nonumber 
\end{align}
together with the analogous expressions obtained by applying $\varpi$ to the elements in eq.~\eqref{d-block}. As described above, the $\tM$-bimodule $N_{\tM}$ is obtained by applying the map $\varpi$ to the $\tH$-coinvariant elements of the right comodule $P\otimes(I_P\cap\tM)$. Thus we see that the bimodule $N_{\tM}$ is generated by the elements
\begin{align}\label{NM}
&\varpi(\sa_\alpha{}^j\sa_\beta{}^l\otimes\sc_i{}^\alpha\sc_k{}^\beta),\quad \varpi(\sa_\alpha{}^j\sb_\beta{}^l\otimes\sc_i{}^\alpha\sc_k{}^\beta),\quad\quad \varpi(\sb_\alpha{}^j\sa_\beta{}^l\otimes\sc_i{}^\alpha\sc_k{}^\beta),\nonumber \\
&\varpi(\sb_\alpha{}^j\sb_\beta{}^l\otimes\sc_i{}^\alpha\sc_k{}^\beta),\quad \varpi(\sa_\alpha{}^j\otimes\sc_i{}^\alpha(\sd_k{}^l-\delta_k^l)),\quad \varpi(\sb_\alpha{}^j\otimes\sc_i{}^\alpha(\sd_k{}^l-\delta_k^l)),
\end{align}
together with all elements obtained by applying $\varpi$ to the generators \eqref{d-block}.

The first condition in eq.~\eqref{r-subs} is just the definition \eqref{bmod1} of $N_M$. To obtain the second condition, we compute first of all that
\begin{align*}
(\delta_R\circ\varpi)(1\otimes\sc_i{}^j\sc_k{}^l) &=\left(\varpi(1\otimes\sc_\mu{}^j\sc_\nu{}^l)\right)\otimes S(u_i{}^\mu)S(u_k{}^\nu), \\ 
(\delta_R\circ\varpi)(1\otimes\sc_i{}^j(\sd_k{}^l-\delta_k^l)) &=\left(\varpi(1\otimes\sc_\mu{}^j(\sd_k{}^l-\delta_k^l))\right)\otimes S(u_i{}^\mu).
\end{align*}
Using this, we find for example that
\begin{align*}
(\delta_R\circ\varpi)(\sa_\alpha{}^j\sa_\beta{}^l\otimes\sc_i{}^\alpha\sc_k{}^\beta) &=\left(\varpi(\sa_\alpha{}^j\sa_\beta{}^l\otimes\sc_\mu{}^\alpha\sc_\nu{}^\beta)\right)\otimes  S(u_i{}^\mu)S(u_k{}^\nu), \\
(\delta_R\circ\varpi)(\sa_\alpha{}^j\otimes\sc_i{}^\alpha(\sd_k{}^l-\delta_k^l)) &=\left(\varpi(\sa_\alpha{}^j\otimes\sc_\mu{}^\alpha(\sd_k{}^l-\delta_k^l))\right)\otimes S(u_i{}^\mu), 
\end{align*}
and so the generators $\varpi(\sa_\alpha{}^j\sa_\beta{}^l\otimes\sc_i{}^\alpha\sc_k{}^\beta)$ and $\varpi(\sa_\alpha{}^j\otimes\sc_i{}^\alpha(\sd_k{}^l-\delta_k^l))$ are stable under the right $H$-coaction. Similarly for the other generators of $N_{\tM}$ appearing in \eqref{NM}. The fact that $\delta_R\circ\varpi$ applied to the generators \eqref{d-block} yields elements in $N_{\tM}\otimes H$ follows immediately from the fact that the ideal $I_{\tH}$ defined in Ex.~\ref{ex:4dopp} is $\textup{Ad}_R$-stable.

The third condition in eq.~\eqref{r-subs} is also verified by direct computation. For example
\begin{align*}
\textup{ver}\left(\varpi(1\otimes\sc_i{}^j\sc_k{}^l)\right)&=\left(S(\sa_\mu{}^l)S(\sa_\nu{}^j) \sc_\alpha{}^\nu\sc_\beta{}^\mu+S(\sc_\mu{}^l)S(\sa_\nu{}^j) \sc_\alpha{}^\nu\sd_\beta{}^\mu  \right. \\
&\qquad\qquad \left. +S(\sa_\mu{}^l)S(\sc_\nu{}^j) \sd_\alpha{}^\nu\sc_\beta{}^\mu + S(\sc_\mu{}^l)S(\sc_\nu{}^j)  \sd_\alpha{}^\nu\sd_\beta{}^\mu\right)\otimes S(u_i{}^\alpha)S(u_k{}^\beta) \\
&=\ep_P(\sc_\alpha{}^j)\ep_P(\sc_\beta{}^k)\otimes S(u_i{}^\alpha)S(u_k{}^\beta)=0,
\end{align*}
with a similar computation yielding that $\textup{ver}\left(\varpi(1\otimes(\sc_i{}^j(\sd_k{}^l-\delta_k^l)))\right)=0$. We deduce that all of the generators \eqref{NM} are zero in the image of the map $\textup{ver}$. 
Finally, for each of the generators $\sd_i{}^j$, $j=1,2$, we have
\begin{align*}
\textup{ver}\circ\varpi(1\otimes\sd_i{}^j)&=\textup{ver}(S(\sb_\mu{}^j)\otimes\sc_i{}^\mu+S(\sd_\mu{}^j)\otimes \sd_i{}^\mu) \\
&= (S(\sb_\mu{}^j)\sc_\alpha{}^\mu+S(\sd_\mu{}^j) \sd_\alpha{}^\mu)\otimes S(u_i{}^\alpha) \\
&=\ep_P(\sd_\alpha{}^j)\otimes S(u_i{}^\alpha)=1\otimes S(u_i{}^j),
\end{align*}
from which it immediately follows that $\textup{ver}(N_{\widetilde M})=\tM\otimes I_H$.
\endproof

Thus we have equipped the quantum four-sphere with a first order differential calculus $\Omega^1S^4_q$ which is compatible with that of the total space of the Hopf fibration $S^7_q\to S^4_q$. The next theorem gives this calculus $\Omega^1S^4_q$ a more concrete description.

\begin{thm}\label{th:fr-sph}
The quantum space $\A[S^4_q]=\A[S^7_q]^H$ has cotangent bundle
\begin{equation}\label{sph-iso}
\Omega^1S^4_q\simeq (P\otimes V)^{{\widetilde H}\otimes H},
\end{equation}
where the right ${\widetilde H}\otimes H$-comodule $V$ is determined to be
\begin{equation}\label{V4s}
V=\textup{Span}_{\C}\,\left\{ x_1,x_1^*,x_2,x_2^* \right\}.
\end{equation}
\end{thm}
\proof 
The algebra $M=\A[S^4_q]$ is generated
by the elements $x_1$, $x_1^*$, $x_2$, $x_2^*$ and $x_0=x_0^*$,
which we identify with generators of the $H$-coinvariant subalgebra
of $\A[S^7_q]$ as in Prop.~\ref{pr:proj}.
Following the usual strategy for a framing comodule, we define the vector space 
$V$ to be the quotient $P^+\cap M/I_P\cap M$. 
Since the generator $x_0$
is not in the kernel of the counit $\ep_P$, the ideal $P^+\cap M$ is
generated by $x_1$, $x_1^*$, $x_2$, $x_2^*$. Since the generators of
$M$ are just linear combinations of the generators of $\A[\CP^3_q]$,
it is easy to check that all quadratic elements in $P^+\cap M$ are
also in $I_P\cap M$, so that $V$ is the
four-dimensional vector space given above.

From eqs.~\eqref{matsum} and \eqref{lin iso} we find that, up to scalar multiples, the elements $x_i$, $x_j^*$, $i,j=1,2$ are given by $v^{kl}:=\iq \sc_1{}^k\sd_2{}^l-\sc_2{}^k\sd_1{}^l$ for  $k,l=1,2$. The vector space $V$ carries a right $\widetilde H$-coaction inherited from that on the vector space $\widetilde V$ in \eqref{V7s}, namely
\begin{equation}\label{Vco1}
V\to V\otimes \widetilde{H},\qquad v^{kl}=(\iq \sc_1{}^k\sd_2{}^l-\sc_2{}^k\sd_1{}^l)~\mapsto~(\iq \sc_1{}^\mu\sd_2{}^l-\sc_2{}^\mu\sd_1{}^l)\otimes \tS (\sa_\mu{}^k),
\end{equation}
for each $k,l=1,2$. We further equip $V$ with a right $H$-coaction by setting
\begin{equation}\label{Vco2}
V\to V\otimes H,\qquad v^{kl}=(\iq \sc_1{}^k\sd_2{}^l-\sc_2{}^k\sd_1{}^l)~\mapsto~(\iq \sc_1{}^k\sd_2{}^\nu-\sc_2{}^k\sd_1{}^\nu)\otimes S^2(u_\nu{}^l),
\end{equation}
for each $k,l=1,2$. Note the use of the square of the antipode in the last formula to ensure that we obtain the correct collection of coinvariant elements below.

These two coactions commuting, we think of $V$ as a right $\widetilde{H}\otimes H$-comodule. The vector space $P=\A[\Sp_q(2)]$ is a right $\widetilde{H}\otimes H$-comodule as well. We equip $P\otimes V$ with the right tensor product coaction. Proving the theorem is thus a matter of checking that the map 
$$\varpi:(P\otimes V)^{\widetilde{H}\otimes H}\to \Omega^1S^4_q
$$ 
is well-defined and a bimodule isomorphism.
The map $\varpi$ is by construction an $M$-$M$ bimodule map. Since the vector space $V$ is spanned by the vectors $v^{kl}$, for $k,l=1,2$, the space $(P\otimes V)^{\widetilde{H}\otimes H}$ is spanned as a left $M$-module by the elements
\begin{equation}\label{els}
\sa_\mu{}^k\sc_\nu{}^l\otimes v^{\mu\nu},\qquad \sa_\mu{}^k\sd_\nu{}^l\otimes v^{\mu\nu},\qquad \sb_\mu{}^k\sc_\nu{}^l\otimes v^{\mu\nu},\qquad \sa_\mu{}^k\sd_\nu{}^l\otimes v^{\mu\nu},
\end{equation}
for each $k,l=1,2$. It is straightforward to check that these are indeed ${\widetilde{H}\otimes H}$-coinvariant, although they are not all independent. Furthermore, by direct computation we find that
\begin{align*}
1\otimes &(\iq \sc_1{}^k\sd_2{}^l-\sc_2{}^k\sd_1{}^l) - (\iq \sc_1{}^k\sd_2{}^l-\sc_2{}^k\sd_1{}^l)\otimes 1= \\
&=\varpi\left( \iq(\sa_\mu{}^k\sb_\nu{}^l\otimes\sc_1{}^\mu\sc_2{}^\nu+\sc_\mu{}^k\sb_\nu{}^l\otimes\sd_1{}^\mu\sc_2{}^\nu+
\sa_\mu{}^k\sd_\nu{}^l\otimes\sc_1{}^\mu\sd_2{}^\nu+\sc_\mu{}^k\sd_\nu{}^l\otimes\sd_1{}^\mu\sd_2{}^\nu) \right. \\
&\quad\quad -\left(\sa_\mu{}^k\sb_\nu{}^l\otimes\sc_2{}^\mu\sc_1{}^\nu+\sc_\mu{}^k\sb_\nu{}^l\otimes\sd_2{}^\mu\sc_1{}^\nu
+\sa_\mu{}^k\sd_\nu{}^l\otimes\sc_2{}^\mu\sd_1{}^\nu+\sc_\mu{}^k\sd_\nu{}^l\otimes\sd_2{}^\mu\sd_1{}^\nu\right) \\
&\quad\quad\quad\left.-(\iq\sc_1{}^k\sd_2{}^l-\sc_2{}^k\sd_1{}^l)\otimes 1)\right) \\
&=\varpi\left(\sa_\mu{}^k\sd_\nu{}^l\otimes (\iq\sc_1{}^\mu\sd_2{}^\nu-\sc_2{}^\mu\sd_1{}^\nu)+\sc_\mu{}^k\sb_\nu{}^l\otimes (\iq\sd_1{}^\mu\sc_2{}^\nu-\sd_2{}^\mu\sc_1{}^\nu) \right. \\
&\quad\quad \left.-\iq\sc_1{}^k\sd_2{}^l\otimes (q\sd_2{}^1\sd_1{}^2-\sd_1{}^1\sd_2{}^2+1)+\sc_2{}^k\sd_1{}^l\otimes (\iq\sd_1{}^2\sd_2{}^1 - \sd_2{}^2\sd_1{}^1+1)\right) \\
&=\varpi\left(\sa_\mu{}^k\sd_\nu{}^l\otimes(\iq\sc_1{}^\mu\sd_2{}^\nu-\sc_2{}^\mu\sd_1{}^\nu)+\sc_\nu{}^k\sb_\nu{}^l\otimes  (\iq\sd_1{}^\mu\sc_2{}^\nu-\sd_2{}^\mu\sc_1{}^\nu) \right. \\
&\quad\quad \left.-\iq\sc_1{}^k\sd_2{}^l\otimes(1-\tfrac{1}{2}(1-x_0))+\sc_2{}^k\sd_1{}^l\otimes(1-\tfrac{1}{2}(1-x_0))\right) \\
&=\varpi\left(\sa_\mu{}^k\sd_\nu{}^l\otimes (\iq\sb_\mu{}^2\sd_\nu{}^1-\sb_\mu{}^1\sd_\nu{}^2) +\sc_\nu{}^k\sb_\nu{}^l\otimes  (\iq\sd_1{}^\mu\sc_2{}^\nu-\sd_2{}^\mu\sc_1{}^\nu) \right. \\
&\quad\quad \left.+(\sc_2{}^k\sd_1{}^l-\iq\sc_1{}^k\sd_2{}^l)\otimes(\sc_2{}^2\sc_1{}^1-\sc_1{}^2\sc_2{}^1)\right) \\
&=\varpi\left(\sa_\mu{}^k\sd_\nu{}^l\otimes (\iq \sc_1{}^\mu\sd_2{}^\nu-\sc_2{}^\mu\sd_1{}^\nu) +\sc_\nu{}^k\sb_\nu{}^l\otimes  (\iq\sd_1{}^\mu\sc_2{}^\nu-\sd_2{}^\mu\sc_1{}^\nu)\right) ,
\end{align*}
for each $k,l=1,2$. We have used the fact that products of the form $\sc_i{}^j\sc_k{}^l$ are equivalent to zero in the quotient space $V$. Using the commutation relations in the quantum group $P$, it follows that certain linear combinations of the elements \eqref{els} are carried by the map $\varpi$ onto the one-forms
$$
\D x_i=1\otimes x_i-x_i\otimes 1,\qquad \D x_j^*=1\otimes x_j^*- x_j^*\otimes 1,\qquad i,j=1,2,
$$
so the map $\varpi$ is surjective onto $\Omega^1S^4_q$. By construction, the kernel of $\varpi$ applied to the bimodule $(P\otimes (P^+\cap M))^{\tH\otimes H}$ is precisely the sub-bimodule $(P\otimes (I_P\cap M))^{\tH\otimes H}$, which is by definition carried onto the sub-bimodule $N_M$ of $\widetilde\Omega^1S^4_q$, which defines the calculus $\Omega^1S^4_q$, so that $\varpi$ is also injective when considered as a map from $(P\otimes V)^{\tH\otimes H}$ onto $\Omega^1S^4_q$.
\endproof

Of course, one could proceed further and calculate the bimodule
relations in the calculus $\Omega^1S^4_q$, as well as those in the
calculi $\Omega^1S^7_q$ and $\Omega^1\CP^3_q$. However, for us the important
facts are that the calculus $\Omega^1S^4_q$ is four-dimensional and
explicitly described by the isomorphism \eqref{sph-iso}. In the next
section, we shall use this knowledge to equip the quantum sphere
with a Hodge structure on two-forms and introduce corresponding 
anti-self-duality equations.

\section{The Instanton Solution}\label{se:bi}
We need to study the differential structure of the quantum sphere
$S^4_q$ in further detail. Having equipped the quantum sphere $S^4_q$ with a first order differential structure, 
we have to extend it to obtain higher order differential forms. We will then be able to formulate the notion of a Hodge
structure on the quantum four-sphere and state an associated set of
anti-self-duality equations.
As we shall see, this is a matter of understanding the representation theory 
 of the quantum group $\SU^\cop_q(2)\times\SU_q(2)$, out of which we construct a suitable operator leading to 
 an `antisymmetric' tensor algebra of forms.

\subsection{Representation theory of $\SU^\cop_q(2)\times\SU_q(2)$}

We continue to write $M=P^{\tH\otimes H}$ for the quantum
principal bundle with 
$$
P=\A[\Sp_q(2)], \qquad \tH=\A[\SU^\cop_q(2)],\qquad H=\A[\SU_q(2)],\qquad M=\A[S^4_q].
$$
We write $V_j$ and $\tV_j$ to denote the irreducible spin $j$ corepresentations of 
$H=\A[\SU_q(2)]$ and $\tH=\A[\SU^\cop_q(2)]$, respectively. In particular, we need the  spaces $V_{\frac{1}{2}}\otimes V_{\frac{1}{2}}$ and $\tV_{\frac{1}{2}}\otimes \tV_{\frac{1}{2}}$, equipped with the respective
right tensor product coactions of $\A[\SU_q(2)]$ and $\A[\SU^\cop_q(2)]$.

\begin{lem}\label{le:br-iso}
There is an isomorphism 
$$
\Psi:V_{\frac{1}{2}}\otimes V_{\frac{1}{2}}\to V_{\frac{1}{2}}\otimes V_{\frac{1}{2}}
$$ 
of right $\A[\SU_q(2)]$-comodules defined by
\begin{align*}
\Psi(x\otimes x)&=q^{1/2}x\otimes x, \\
 \Psi(x\otimes y)&=q^{-1/2}\left( y\otimes x+(q-\iq)x\otimes y\right), \\
\Psi(y\otimes x)&=q^{-1/2}x\otimes y, \\
\Psi(y\otimes y)&=q^{1/2}y\otimes y,
\end{align*}
whose eigenvalues are
$q^{1/2}$ and $-q^{-3/2}$ with multiplicities $3$ and $1$
respectively.
\end{lem}

\proof The given map is nothing other than the braiding on the category of right $\ASU$-comodules induced by the coquasitriangular structure \eqref{ct}, which means precisely that it is an isomorphism as claimed \cite{Ma:book}. If in doubt, a direct verification is not difficult. 
The matrix $R$ occuring in \eqref{ct} obeys the Hecke relation
$(R-q)(R+\iq)=0$, from which its eigenvalues are easily found to
be $q$ and $-\iq$ with multiplicities $3$ and $1$ respectively. 
Due to the normalisation factor $\zeta=q^{-1/2}$, the eigenvalues of
$\Psi$ are therefore just $q^{1/2}$ and $-q^{-3/2}$ with the stated multiplicities.
\endproof

Similarly, there is an isomorphism
$$
\widetilde{\Psi}:\tV_{\frac{1}{2}}\otimes \tV_{\frac{1}{2}}\to \tV_{\frac{1}{2}}\otimes \tV_{\frac{1}{2}}
$$
of right $\A[\SU^\cop_q(2)]$-comodules, determined by the canonical braiding on the category of such comodules induced again by the coquasitriangular structure \eqref{c-rels}. The eigenvalues here are once again $q^{1/2}$ and $-q^{-3/2}$ with multiplicities $3$ and $1$ respectively.

Writing $W:=\tV_{\frac{1}{2}}\otimes V_{\frac{1}{2}}$, it follows that we have an isomorphism
$$
\Phi:W\otimes W\to W\otimes W
$$
of right $\tH\otimes H$-comodules defined by the formula
$$
\Phi:=(\id\otimes \sigma \otimes
\id)\circ(\widetilde\Psi\otimes\Psi)\circ(\id\otimes \sigma \otimes \id),
$$
where $\sigma:V_{\frac{1}{2}}\otimes \tV_{\frac{1}{2}}\to \tV_{\frac{1}{2}}\otimes V_{\frac{1}{2}}$
denotes the ordinary `flip' on tensor factors. By inspection we see that the eigenvalues of the map $\Phi$ are $q$, $-\iq$ and $q^{-3}$ with multiplicities
$9$, $6$ and $1$ respectively. We would like to use the operator $\Phi$ to define the notion of  an `antisymmetric' tensor algebra over $W$.

\begin{defn}\label{2d-anti}
We say that an element
$v\in W\otimes W$ is {\em $q$-antisymmetric} if it belongs to
the six-dimensional eigenspace of $\Phi$ of the eigenvalue $-\iq$.
\end{defn}

To generalize this to elements of higher order, we need some additional manipulations. It is straightforward to verify that $\Phi$ obeys the
braid relation
\begin{equation}\label{braid}
(\id\otimes\Phi)(\Phi\otimes\id)(\id\otimes\Phi)=(\Phi\otimes\id)(\id\otimes\Phi)(\Phi\otimes\id)
\end{equation}
on $W\otimes W\otimes W$, as one might expect. More generally, for each $r=2,3,\ldots$, we have a set
\begin{equation}\label{br-set}
\{\Phi_1,\Phi_2,\ldots,\Phi_{r-1}\}
\end{equation}
of automorphisms of the $\tH\otimes H$-comodule $\otimes^r\,W:=W\otimes
\cdots\otimes W$ made of $r$ copies of $W$ and equipped with
the tensor product coaction. The automorphism $\Phi_k$ is
defined by 
$$
\Phi_k:=\id\otimes\cdots\otimes
\Phi\otimes\cdots\otimes\id,
$$
where the tensor product has $r-1$ factors and the map $\Phi$ occurs
in the $k$-th position. It is immediate from eq.~\eqref{braid} that
these $\Phi_k$ satisfy the braid relation
\begin{equation}\label{br-rel}
\Phi_k\Phi_{k+1}\Phi_k=\Phi_{k+1}\Phi_k\Phi_{k+1}, \qquad \textup{for} \qquad k=1,2\ldots,r-1 .
\end{equation}

Following \cite{wor}, let $\Sigma_r$ denote the permutation group on
$r$ objects and consider the set of nearest neighbour transpositions 
\begin{equation}\label{trans}
\{t_1,t_2,\ldots,t_{r-1}\},
\end{equation}
where the permutation $t_k$ is defined for each $k=1,2,\ldots,r-1$
as the operation which exchanges the object in position $k$ with the one in position $k+1$, whilst leaving all
other objects in their places. For each $p\in\Sigma_r$ we denote by
$I(p)$ the number of pairs of elements in the sequence
$(p(1),p(2),\ldots,p(r))$ for which $i<j$ but $p(j)<p(i)$; then the permutation $p$ is a product of
$I(p)$ elements in the set \eqref{trans},
\begin{equation}\label{decomp}
p=t_{k_1}t_{k_2}\ldots t_{k_{I(p)}}.
\end{equation}
We define an automorphism $\Pi_p$ of the $\tH\otimes H$-comodule $\otimes^r W$ by
replacing each of the transpositions in the permutation
\eqref{decomp} by the corresponding element of the set
\eqref{br-set}. Due to the braid relation \eqref{br-rel}, the
resulting automorphism is independent of the decomposition
\eqref{decomp}. We have then,
\begin{equation}\label{br-auto}
\Pi_p:=\Phi_{k_1}\Phi_{k_2}\ldots
\Phi_{k_{I(p)}}.
\end{equation}
Using these, finally we define the `antisymmetrisation operator' $A_r$ by the
formula
\begin{equation}\label{antisym}
A_r:=\sum_{p\in\Sigma_r} (-1)^{I(p)}\Pi_p.
\end{equation}
For $r=1$ we set $A_1:=\id$, the identity operator on $W$. 
Also, we set $W^0:=\C$ to be the trivial $\tH\otimes H$-comodule.
This means that we arrive at the following definition.

\begin{defn}\label{q-forms}
The vector space $W^{\wedge r}$ of {\em $q$-antisymmetric} elements in the tensor product $\otimes^r W$ is defined to be
$$
W^{\wedge r}:=\otimes^r W/\mathrm{Ker}\,A_r, \qquad r=1,2,\ldots.
$$
Moreover, for each $v\in W^{\wedge r}$ and $w\in W^{\wedge s}$ we define $v\wedge_q w \in 
W^{\wedge (r+s) }$ to
be the projection of the tensor $v\otimes w$ to its
$q$-antisymmetric part.

\end{defn}

A close inspection finds the vector spaces
$W^{\wedge r}$ to be of dimensions $1$, $4$, $6$, $4$ and $1$ when $r=0,1,2,3$ and $4$,
respectively, and zero otherwise.
We write $W^\bullet:=\oplus_r W^{\wedge r}$ for the graded vector space of all $q$-antisymmetric
tensors. Since the
$\tH\otimes H$-coaction preserves the eigenspaces of $\Phi$, the
tensor product coaction descends to a right
$\tH\otimes H$-coaction on each of the quotient vector spaces $W^{\wedge r}$ and hence to a coaction on $W^\bullet$.

\begin{lem}\label{le:assoc}
The wedge product $v\otimes w\mapsto v\wedge_q w$ is $\tH\otimes H$-covariant
and makes $W^\bullet$ into a graded associative algebra.
\end{lem}

\proof The fact that the product is graded and associative is
straightforward. Covariance of the product is precisely the statement that
the braiding $\Phi$ is an intertwiner for the $\tH\otimes H$-coaction
on the tensor product $W=\tV_{\frac{1}{2}}\otimes V_{\frac{1}{2}}$ and so preserves its
eigenspaces.\endproof

\subsection{Hodge structure on $S^4_q$}\label{se:hodge} Using the combinatorics of the previous section, we are ready to equip $S^4_q$ with a full algebra of differential forms.
To this end, again recall from Theorem~\ref{th:fr-sph} that the sphere $S^4_q$ is a 
quantum manifold having cotangent bundle
\begin{equation}\label{s-frame}
\Omega^1S^4_q\simeq (P\otimes V)^{\tH\otimes H},
\end{equation}
with $V$ given in \eqref{V4s} and corresponding and rather complicated right $\tH\otimes H$-coaction in eqs. \eqref{Vco1} and 
\eqref{Vco2}.
The following lemma makes things a little easier to work with.

\begin{lem}\label{le:V-iso}
There is a unitary equivalence
\begin{equation}\label{V}
V\simeq W:=\tV_{\frac{1}{2}}\otimes V_{\frac{1}{2}}
\end{equation}
of right $\A[\SU_q^\cop(2)]\otimes\A[\SU_q(2)]$-comodules.
\end{lem}

\proof The vector space $V$ in \eqref{V4s} is a four-dimensional irreducible right
$\tH\otimes H$-comodule. Up to unitary equivalence, the irreducible
four-dimensional $\tH\otimes H$-comodules are just 
$$\tV_{\frac{3}{2}}\otimes V_0,\qquad \tV_{\frac{1}{2}}\otimes V_{\frac{1}{2}},\qquad \tV_{0}\otimes V_{\frac{3}{2}}.
$$ 
By inspection we see that both of the tensor factors $\A[\SU_q^\cop(2)]$
and $\A[\SU_q(2)]$ coact non-trivially on $V$, whence we must have the
isomorphism as stated.\endproof

In the notation of eq.~\eqref{assoc-bd} for associated vector bundles, 
the cotangent bundle $\Omega^1S^4_q$ is isomorphic to the vector bundle $\mathcal{M}(V)$ associated to the 
$\tH\otimes H$ principal bundle $\A[S^4_q]\hookrightarrow\A[\Sp_q(2)]$. In the spirit of associated bundles, one is led immediately to the following definition of the differential forms of higher order.

\begin{defn}
We define the space of differential $r$-forms on the quantum sphere
$S^4_q$ to be the $\A[S^4_q]$-bimodule $\Omega^rS^4_q:=\mathcal{M}(V^{\wedge r})$.
\end{defn}

Here $V^{\wedge r}$ is the space of {\em $q$-antisymmetric} elements as in Definition~\ref{q-forms}.
Next we make $\Omega^\bullet S^4_q:=\oplus_r \Omega^r S^4_q$ into a
graded associative algebra. Indeed, the wedge product $\wedge_q$ on $V^\bullet$ as in Lemma~\ref{le:assoc},
makes it easy to equip
$\Omega^\bullet S^4_q$ with a graded associative algebra structure.
For $\omega=a\otimes v$ and $\omega'=b\otimes w$, with $a,b\in
P=\A[\Sp_q(2)]$ and $v\in V^{\wedge r}$, $w\in V^{\wedge s}$, we define
$$
\omega\wedge_q\omega':=(ab)\otimes(v\wedge_q w).
$$
This gives a well-defined product on differential
forms.
\begin{lem}\label{le:d-assoc}
The product $\wedge_q$ makes $\Omega^\bullet S^4_q$ into a graded
associative algebra.
\end{lem}

\proof 
Associativity of the product on $\Omega^\bullet S^4_q$
follows from associativity of the product on $V^\bullet$.
It remains to check that, if $a,b\in P=\A[\Sp_q(2)]$ and $v\in
V^{\wedge r}$, $w\in V^{\wedge s}$ are such that $a\otimes v$ and $b\otimes w$ are
$\tH\otimes H$-coinvariant, then so is the product $(ab)\otimes(v\wedge_q w)$.
To this end, we compute that under the $\tH\otimes H$-coaction 
\begin{align*}
(ab)\otimes(v\wedge_q w)&\mapsto
(ab)\bz\otimes(v\wedge_q w)\bz\otimes (ab)\uo(v\wedge_q w)\uo \\
& \qquad = a\bz b\bz\otimes (v\bz\wedge_q w\bz)\otimes (a\uo v\uo)(b\uo w\uo)\\
& \qquad =(ab)\otimes (v\wedge_q w)\otimes 1,
\end{align*}
as required. 
\endproof

The product $\wedge_q$ is surjective (which may be checked just as in \cite{l-da}, for example) and so the calculus
$\Omega^\bullet S^4_q$ is generated in degree one, meaning that every
form of degree $r\geq 1$ can be written as a product of one-forms.
From general considerations \cite{wor}, it is automatic that the exterior derivative
$\D:\A[S^4_q]\to\Omega^1S^4_q$ (which we have not explicitly given, since we do not need it) extends uniquely to a de Rham complex
$\D:\Omega^rS^4_q\to\Omega^{r+1}S^4_q$ by requiring it to satisfy
$\D^2=0$ and the graded Leibniz rule
\begin{equation}\label{gr-leib}
\D(\omega\wedge_q\omega')=(\D\omega)\wedge_q\omega'+(-1)^r\omega\wedge_q(\D\omega'),\qquad
\omega\in\Omega^rS^4_q,~\omega'\in\Omega^{r'}S^4_q.
\end{equation}
The $*$-structure on $\A[S^4_q]$ also extends uniquely to an
involution on $\Omega^\bullet S^4_q$ satisfying 
$$\D\omega^*=-(\D\omega)^* , \qquad \textup{for all} \qquad \omega\in \Omega^\bullet S^4_q .
$$ 

We are interested, in particular, in the space $\Omega^2S^4_q$ of
two-forms. Just as in the classical case, this space has a nice
decomposition into components given by irreducible representations of $\SU^\cop_q(2)\times\SU_q(2)$.

\begin{prop}\label{pr:sd-asd}
There is a decomposition of $\A[S^4_q]$-bimodules
$$
\Omega^2S^4_q\simeq \Omega^2_+ S^4_q \oplus \Omega^2_- S^4_q,
$$
where the sub-bimodules $\Omega^2_\pm S^4_q$ are to be
determined.
\end{prop}

\proof Just as in the classical case, the required decomposition of the bimodule
$\Omega^2S^4_q$ will come from
decomposing the $\tH\otimes H$-comodule $V^{\wedge 2}$ into irreducible sub-comodules. The
latter is just the space of $q$-antisymmetric tensors in the
comodule
\begin{align*}
V\otimes V&\simeq (\tV_{\frac{1}{2}}\otimes V_{\frac{1}{2}})\otimes
(\tV_{\frac{1}{2}}\otimes V_{\frac{1}{2}})\\
&=(\tV_0\otimes V_0)\oplus(\tV_0\otimes V_1)\oplus(\tV_1\otimes
V_0)\oplus (\tV_1\otimes V_1),
\end{align*}
where we have expanded the tensor products as a Clebsch-Gordan
series of irreducible sub-comodules. It is clear from Definition~\ref{2d-anti} that the
antisymmetric tensors in this case are precisely those in the
subspace $(\tV_0\otimes V_1)\oplus (\tV_1 \otimes V_0)$. We define
$$
\Omega^2_+ S^4_q:=\mathcal{M}(\tV_1\otimes V_0)\simeq \mathcal{M}(\tV_1),\qquad
\Omega^2_- S^4_q:=\mathcal{M}(\tV_0\otimes V_1)\simeq \mathcal{M}(V_1)
$$
to obtain the decomposition as stated.\endproof

\begin{defn}
The {\em Hodge $*$-operator} $*:\Omega^2S^4_q\to \Omega^2S^4_q$ is
the linear map defined by 
$$
*(\omega_\pm):=\pm\omega_\pm \qquad \textup{for} \quad \omega_\pm\in\Omega^2_\pm S^4_q ,
$$  
and extended by left
$\A[S^4_q]$-linearity. A two-form $\omega_+\in\Omega^2_+ S^4_q$ is said to
be {\em self-dual}; a two-form $\omega_-\in\Omega^2_- S^4_q$ is said to be
{\em anti-self-dual}.
\end{defn}

\subsection{Differential structure of twistor space}
Having studied the higher-order differential forms on the sphere
$S^4_q$, we turn to the differential geometry of twistor space
$\CP^3_q$ that we construct in a similar manner. 
This time we start with the framing of twistor space
$\CP^3_q$ computed in Theorem~\ref{th:tw-frm}, with the quantum 
principal bundle $\hM=P^{H'\otimes H}$ now given by  
$$
P=\A[\Sp_q(2)], \qquad H'=\A[\U(1)], \qquad H=\A[\SU_q(2)],\qquad\hM=\A[\CP^3_q]
$$
and corresponding cotangent bundle $\Omega^1\CP^3_q\simeq (P\otimes\widehat V)^{\hH}$.
In the notation of \S\ref{se:su2}, we continue to write $V_j$ for the irreducible $\A[\SU_q(2)]$-comodule with spin $j$ while $U_k$ denotes the irreducible $\A[\U(1)]$-comodule with index $k$.

\begin{prop}
There is a unitary equivalence of right $H'\otimes H$-comodules 
\begin{equation}\label{V-eq}
\hV\simeq (U_2\otimes V_0)\oplus(U_1\otimes V_{\frac{1}{2}})\oplus(U_0\otimes
V_0)\oplus(U_{-1}\otimes V_{\frac{1}{2}})\oplus(U_{-2}\otimes V_0).
\end{equation}
\end{prop}

\proof This is immediate by inspection of the formul{\ae} \eqref{tw-fr}. Indeed, we have that
\begin{align*}
&U_2\otimes V_0=\textup{Sp}_{\C}\{\qp_{12}\},\qquad U_0\otimes V_0=\textup{Sp}_{\C}\{\qp_{11}\},\qquad U_{-2}\otimes V_0=\textup{Sp}_{\C}\{\qp_{21}\}, \\
&U_1\otimes V_{\frac{1}{2}}=\textup{Sp}_{\C}\{\qp_{32},\qp_{42}\},\qquad U_{-1}\otimes V_{\frac{1}{2}}=\textup{Sp}_{\C}\{\qp_{23},\qp_{24}\},
\end{align*}
from which the result follows.\endproof

We shall also use the shorthand notation
$$
T_+:=(U_2\otimes V_0)\oplus(U_1\otimes V_{\frac{1}{2}}), \quad
T_0:=(U_0\otimes V_0), \quad T_-:=(U_{-1}\otimes
V_{\frac{1}{2}})\oplus(U_{-2}\otimes V_0), 
$$
so that $\hV\simeq T_+\oplus T_0\oplus T_-$. In terms of the decomposition \eqref{tw-2} for $\Omega^1\CP^3_q$, 
clearly we have
\begin{equation}\label{sub}
\Omega^1_+\CP^3_q \simeq\mathcal{M}(T_+),\qquad \Omega^1_0\CP^3_q \simeq\mathcal{M}(T_0),
\qquad \Omega^1_-\CP^3_q\simeq\mathcal{M}(T_-).
\end{equation}

As we did for the four-sphere, we use the associated bundle
construction to define higher-order differential forms on twistor
space.

\begin{lem}
There is an isomorphism $\widehat\Phi:\hV\otimes \hV\to
\hV\otimes \hV$ of right $H'\otimes H$-comodules whose eigenvalues
are $1$, $-1$, $q^{1/2}$ and $-q^{-3/2}$.
\end{lem}

\proof 
Since the Hopf algebra $H'=\A[\U(1)]$ is commutative, the braiding on the set of comodules $U_k$, for $k\in\ZZ$, is just the tensor flip $\sigma:U_k\otimes U_\ell\to U_\ell\otimes U_k$. Similarly, it is clear that the tensor flip $\sigma:V_0\otimes V_{\frac{1}{2}}\to V_{\frac{1}{2}}\otimes V_0$ is an isomorphism of $\A[\SU_q(2)]$-comodules.  Moreover, we know from Lemma~\ref{le:br-iso} that there is an isomorphism $\Psi:V_{\frac{1}{2}}\otimes V_{\frac{1}{2}}\to V_{\frac{1}{2}}\otimes V_{\frac{1}{2}}$ of $H$-comodules.
With respect to the decomposition \eqref{V-eq} of the $H'\otimes H$-comodule $\hV$, the required braiding map is therefore just the tensor flip 
$$
(U_j\otimes V_0)\otimes(U_k\otimes V_0)\to(U_k\otimes V_0)\otimes (U_j\otimes V_0),
$$
together with the map
\begin{align*}
\widehat\Phi:(U_j\otimes V_{\frac{1}{2}})\otimes (U_k\otimes V_{\frac{1}{2}})\to
(U_k\otimes V_{\frac{1}{2}})\otimes (U_j\otimes V_{\frac{1}{2}}),
\end{align*}
defined by
\begin{align*}
\widehat\Phi:=(\id\otimes \sigma \otimes
\id)\circ(\sigma\otimes\Psi)\circ(\id\otimes \sigma \otimes \id).
\end{align*}
The single eigenvalue of the tensor flip $\sigma$ acting on tensor products of comodules $U_k$ is just $1$ (because it only ever acts on one-dimensional vector spaces). On the other hand, the flip operator $\sigma$ on products of the form $V_0\otimes V_{\frac{1}{2}}$ has eigenvalues $\pm 1$. Since the eigenvalues of $\Psi$ are just $q^{1/2}$ and $-q^{-3/2}$, the eigenvalues of $\widehat\Phi$ are those as stated. 
\endproof

In this case we say that an element $v\in \hV\otimes \hV$ is
$q$-antisymmetric if it lies in the direct sum of the $-1$ and
$-q^{-3/2}$ eigenspaces, which is a $21$-dimensional subspace of
$\hV\otimes \hV$.
More generally, we use the braiding operator $\widehat\Phi$ to construct a $q$-antisymmetrisation operator on each of the tensor products $\otimes^r\,\hV$, 
just as we did for the quantum four-sphere: the construction is identical and so we shall not repeat the details.

We write $\hV^{\wedge r}$ for the quotient of the space $\otimes^r\,\hV$ by the kernel of the $q$-antisymmetrisation operator built from the braiding $\widehat\Phi$, and we set $\hV^0:=\C$ to be the trivial $\tH\otimes H$-comodule.
Again, just as was the case for the four-sphere, there is an $\tH\otimes H$-covariant
graded associative product on the vector space
$\hV^\bullet:=\oplus_r \hV^{\wedge r}$ defined for $v\in \hV^{\wedge r}$, $w\in \hV^{\wedge s}$
by $v\otimes w\mapsto v\wedge_q w$, where $v\wedge_q w$ is the
$q$-antisymmetric part of $v\otimes w$. 
\begin{defn}
The space of differential $r$-forms on twistor space
$\CP^3_q$ is the $\A[\CP^3_q]$-bimodule
$\Omega^r\CP^3_q:=\mathcal{M}(\hV^{\wedge r})$.
The vector space $\Omega^\bullet\CP^3_q:=\oplus_r\Omega^r\CP^3_q$ is a graded associative algebra, with product defined
by
$$
(a\otimes v)(b\otimes w)=(ab)\otimes(v\wedge_q w),\qquad 
v\in \hV^{\wedge r}, ~w\in \hV^{\wedge s}, ~~a,b\in\A[\Sp_q(2)].
$$
\end{defn}

In exactly the same way as in Lemmata~\ref{le:assoc} and
\ref{le:d-assoc}, one checks that this product is
associative and well-defined on differential forms, besides being $\tH\otimes H$-covariant. 
Moreover, the
exterior derivative $\D:\A[\CP^3_q]\to\Omega^1\CP^3_q$ (which again we have not written explicitly) extends
uniquely to a de Rham complex
$\D:\Omega^r\CP^3_q\to\Omega^{r+1}\CP^3_q$ such that $\D^2=0$ and
$$
\D(\omega\wedge_q\omega')=(\D\omega)\wedge_q\omega'+(-1)^r\omega\wedge_q(\D\omega'),\qquad
\omega\in\Omega^r\CP^3_q, ~~\omega'\in\Omega^{r'}\CP^3_q.
$$
The $*$-structure on $\A[\CP^3_q]$ also extends uniquely to an
involution on $\Omega^\bullet \CP^3_q$ satisfying $\D
\omega^*=-(\D\omega)^*$ for all $\omega\in \Omega^\bullet \CP^3_q$.
The product on $\hV^\bullet$ being surjective means that so is the product
on $\Omega^\bullet\CP^3_q$ and every differential $r$-form can be
expressed as the product of exactly $r$ differential one-forms. Furthermore,
the decomposition 
$$
\Omega^1\CP^3_q=\Omega^1_+\CP^3_q \oplus\Omega^1_0\CP^3_q \oplus\Omega^1_-\CP^3_q
$$ 
in terms of the differential sub-calculi \eqref{sub} gives some extra structure to this differential calculus in the form
of a $\ZZ^3$-grading, defined in the following way.

\begin{defn}
With $(m,a,n)\in \ZZ^3$, an $r$-form $\omega\in\Omega^p\CP^3_q$ is said to be of type $(m,a,n)$ if it is the
product of $m$ elements of $\Omega^1_+\CP^3_q$, $a$ elements of $\Omega^1_0\CP^3_q$ and $n$ elements of 
$\Omega^1_-\CP^3_q$.
\end{defn}

It is clear that an $r$-form of type $(m,a,n)$ must obey $r=m+a+n$.
Moreover, although the values $(m,a,n)$ are allowed {\em a priori}
to take arbitrary value, it is evident that the only non-zero spaces
of forms of type $(m,a,n)$ are those with $m,n\in\{0,1,2,3\}$ and
$a\in \{0,1\}$.

\subsection{The geometry of the twistor fibration}\label{se:ast}
Recall the $q$-deformed version of the Penrose twistor
fibration $\CP^3\to S^4$ expressed as an inclusion of coordinate algebras
$$
\eta:\A[S^4_q]\hookrightarrow\A[\CP^3_q]
$$
in eq.~\eqref{tw-bund}.  
In the classical case, differential forms
on the base space $S^4$ may be pulled back along the fibration to
give differential forms on twistor space $\CP^3$. In the noncommutative case, there is a similar phenomenon defined as follows. 

\begin{prop}
There is an injective linear `push-out' map
$$
\eta_*:\Omega^\bullet S^4_q\to \Omega^\bullet \CP^3_q
$$
which is an intertwiner for the respective $\A[S^4_q]$- and $\A[\CP^3_q]$-bimodule structures.
\end{prop}

\proof Using the canonical projection $\pi:\A[\SU^\cop_q(2)]\to\A[\U(1)]$, the $\A[\SU_q^\cop(2)]$-comodule $\tV_{\frac{1}{2}}$ is automatically an $\A[\U(1)]$-comodule according to
$$
\Delta_{\frac{1}{2}}':\tV_{\frac{1}{2}}\to \tV_{\frac{1}{2}}\otimes\A[\U(1)],\qquad 
\Delta_{\frac{1}{2}}':=(\id\otimes\pi)\circ\widetilde{\Delta}_{\frac{1}{2}}.
$$
By inspection, we see that there is a simple isomorphism $\tV_{\frac{1}{2}}\simeq U_1\oplus U_{-1}$ as $\A[\U(1)]$-comodules. Consequently, there is an injective morphism of right comodules
$$
j:V\to \hV,\qquad V=\tV_{\frac{1}{2}}\otimes V_{\frac{1}{2}}\to (U_1\otimes V_{\frac{1}{2}})\oplus(U_{-1}\otimes V_{\frac{1}{2}})\subseteq \hV,
$$
so that the map $j:V\to \hV$ simply takes the $\A[\SU^\cop_q(2)]\otimes \A[\SU_q(2)]$-comodule $V$ and views it as an $\A[\U(1)]\otimes\A[\SU_q(2)]$-sub-comodule of $\hV$ according to the decomposition \eqref{V-eq}. The operation $V\to V^{\wedge \bullet}$ being functorial there is an induced comodule morphism
$$
j_*: V^{\wedge r} \to \hV^{\wedge r}, \qquad r=0,1,2,\ldots \, .
$$
This induces a similar morphism $\eta_*: P\otimes V^{\wedge r}\to P\otimes\hV^{\wedge r}$ of tensor product comodules and
there is an induced map $\eta_*:\Omega^r S^4_q\to\Omega^r\CP^3$ upon passing to the coinvariant elements.\endproof 

The above simple result gives rise to the following characterisation of anti-self-dual two forms on $S^4_q$. Save for the fact that we are working with $q$-deformed tensors, its proof is just representation theory and is more or less identical to the classical situation ({\em cf}. \cite{ma:gymf}).

\begin{prop}\label{pr:asd}
A two-form $\omega\in\Omega^2S^4_q$ is anti-self-dual if and only if
its push-out $\eta_*\omega \in \Omega^2 \CP^3_q$ along the twistor fibration is a
two-form of type $(1,0,1)$.
\end{prop}

\proof Recall from Prop.~\ref{pr:sd-asd} that the space
$\Omega^2_- S^4_q$ of anti-self-dual two-forms on $S^4_q$ is the space of
two-forms $\mathcal{M}(V_1)$ associated to the spin $1$ corepresentation of the Hopf algebra $\A[\SU_q(2)]$. To prove the claim, we need to examine the image
of the space $\Omega^2_- S^4_q$ under the push-out map $\eta_*$, which we do
by looking at the decomposition of the vector space $\hV$ into
irreducible $\A[\SU_q(2)]$-comodules. Clearly the push-out of a
two-form on $S^4_q$ cannot involve the `extra dimension'
$\mathcal{M}(T^0)$ in the calculus $\Omega^\bullet\CP^3_q$, so the
push-out $\eta_*\omega$ of an anti-self-dual two-form
$\omega\in\Omega^2_- S^4_q$ must have a decomposition
$$
\eta_*\omega=\omega^{2,0,0}+\omega^{1,0,1}+\omega^{0,0,2}
$$
by type. The spaces of two-forms of type $(2,0,0)$ and $(0,0,2)$ are
both one-dimensional trivial $\A[\SU_q(2)]$-comodules, whence we
must have $\omega^{2,0,0}=\omega^{0,0,2}=0$ and so the image
$\eta_*\Omega^2_- S^4_q$ consists of two-forms of type $(1,0,1)$. The converse is now immediate.\endproof

\subsection{Twistor geometry of the instanton}
We are ready for the definition of an instanton on the quantum
four-sphere $S^4_q$. Firstly, a {\em connection} on a finitely
generated projective right $\A[S^4_q]$-module $\E$ is a linear map
$\n:\E\to\E\otimes_{\A[S^4_q]}\Omega^1S^4_q$ obeying the Leibniz rule
$$
\n(\xi a)=(\n \xi)a + \xi\otimes \D a,\qquad \xi\in \E, ~~a \in\A[S^4_q].
$$
Any such $\n$ extends uniquely to
$$
\n:\E\otimes_{\A[S^4_q]}\Omega^rS^4_q\to \E\otimes_{\A[S^4_q]}
\Omega^{r+1}S^4_q,
$$
defined for all $r\geq 0$ and obeying the graded Leibniz rule
\begin{equation*}
\n(\xi\omega)=(\n\xi)\wedge_q\omega+(-1)^r\xi\otimes(\D \omega)
\end{equation*}
for all $\xi\in\E\otimes_{\A[S^4_q]}\Omega^rS^4_q$ and $\omega\in
\Omega^\bullet S^4_q$. The {\em curvature} of the connection $\n$ is the map
$$
F:=\n^2:\E\to \E\otimes_{\A[S^4_q]}\Omega^{2}S^4_q.
$$
This map $F$ is $\A[S^4_q]$-linear, whence we may think of it as a
two-form on $S^4_q$ taking values in the algebra
$\mathrm{End}_{\A[S^4_q]}(\E)$ of $\A[S^4_q]$-linear endomorphisms
of $\E$ or, more precisely, as an element $F\in \mathrm{Hom}_{\A[S^4_q]}(\E, \E\otimes_{\A[S^4_q]}\Omega^{2}S^4_q)\simeq  \mathrm{End}_{\A[S^4_q]}(\E)\otimes_{\A[S^4_q]}\Omega^{2}S^4_q$.

\begin{defn}\label{de:hso}
A connection $\n:\E\to\E\otimes_{\A[S^4_q]}\Omega^1S^4_q$ is said to
be an {\em instanton} if its curvature $F$ is an anti-self-dual two-form for the Hodge $*$-operator 
$$
\id \otimes * \simeq * : \Omega^2 S^4_q\to\Omega^2S^4_q .
$$
\end{defn}

Associated to the projection $\pp$ of eq.~\eqref{defproj} which
defines the four-sphere $S^4_q$ there is a canonical example of a
noncommutative vector bundle over $S^4_q$. 

\begin{defn}
The {\em canonical instanton vector bundle} on the quantum four-sphere $S^4_q$ is the right $\A[S^4_q]$-module $\E:=\Pp\A[S^4_q]^4$ defined by the complementary projection $\Pp:=\id-\pp$.
\end{defn}

As usual, the vector bundle $\E$ comes equipped with
the canonical Grassmann connection $\n:=\Pp\circ \D$. We claim that
this connection has anti-self-dual curvature.
To prove this, we make use of the geometry of the twistor fibration over $S^4_q$ as described in \S\ref{se:ast}.
Indeed, using Prop.~\ref{pr:asd} we immediately deduce that, with the projection
$\Pp:=\id-\pp$, the curvature of the connection $\n=\Pp\circ \D$ on
$S^4_q$ is an anti-self-dual two-form and hence that this connection
is an instanton, as claimed.

\begin{thm}
The curvature $F=\Pp(\D \Pp)^2$ of the connection $\n=\Pp\circ\D$ is
anti-self-dual, that is to say 
$$
*F=-F ,
$$ 
where $*:\Omega^2 S^4_q\to\Omega^2S^4_q$ is the Hodge $*$-operator as in Definition~\eqref{de:hso}.
\end{thm}

\proof Recall from eq.~\eqref{pr-dec} that the projection $\pp$ decomposes as
$\pp=\sfu\sfu^*$, where the partial isometry $\sfu$ obeys
$\sfu^*\sfu=\mathbbm{1}_2$. Using this, one finds that
\begin{align*}
\D(\sfu \sfu^*)=\Pp(\D\sfu)\sfu^*+\sfu(\D\sfu^*)\Pp,
\end{align*}
and hence in turn that
\begin{align*}
(\D\Pp)\wedge_q(\D\Pp)={\Pp}(\D{\sfu})\wedge_{q}(\D{\sfu}^*)\Pp+{\sfu}(\D{\sfu}^*){\Pp}\wedge_q {\Pp}(\D{\sfu}){\sfu}^*,
\end{align*}
where we have used $\sfu^*\Pp=0=\Pp\sfu$. The second
term in the above expression is identically zero when acting on any
element in the image $\E$ of $\Pp$, whence the curvature $F$ of the
connection $\n$ works out to be
\begin{equation*}
F=\Pp(\D\Pp)^2=\Pp(\D{\sfu})\wedge_q(\D{\sfu}^*)\Pp
\end{equation*}
or, in terms of explicit matrix components,
\begin{equation*}
F_{ab}=\sum_{j,l}~\Pp_{aj}\left((\D z_j)\wedge_q (\D z_l^*)+(\D
Jz_j)\wedge_q (\D J z_l^*)\right)\Pp_{lb}
\end{equation*}
for $a,b=1,\ldots,4$. It is clear by inspection that on twistor
space $\A[\CP^3_q]$ this is a two-form of type $(1,0,1)$, whence by
Prop.~\ref{pr:asd} the curvature $F$ is anti-self-dual, as
claimed.\endproof

\end{document}